\documentclass[11pt, titlepage]{article}
\usepackage{amsmath}    
\usepackage{appendix} 
\usepackage{epsfig}
\usepackage{makeidx}               
\usepackage{array}              
\usepackage{latexsym}
\usepackage{multirow}
\usepackage{comment}
\usepackage{booktabs} 
\usepackage{amssymb}
\usepackage{algorithm} 
\usepackage[noend]{algpseudocode}
\usepackage{graphicx}
\usepackage{xcolor} 
\usepackage{pgfgantt}

\graphicspath{ {./images/} } 
\usepackage{color,soul}
\usepackage{amssymb}

\usepackage{natbib}
\bibpunct[, ]{(}{)}{,}{a}{}{,}%
\def\bibfont{\small}%
\def\bibsep{\smallskipamount}%
\def\bibhang{24pt}%

\allowdisplaybreaks

\usepackage{tikz,pgfgantt,setspace}
\newtheorem{proposition}{Proposition}
\newcommand*{\QEDB}{\hfill\ensuremath{\square}}%

\begin{document}
\title{A Stochastic Programming Approach to Surgery Scheduling under Parallel Processing Principle}
\date\today
\begin{titlepage}
	\vspace*{18ex} {\centerline{\bf{A STOCHASTIC PROGRAMMING APPROACH TO
	}}} {\centerline{\bf{SURGERY SCHEDULING UNDER PARALLEL PROCESSING PRINCIPLE}}}\vspace{0.3in}

	\begin{center} {\bf Batuhan Celik$^1$, Serhat Gul$^2$, Melih Celik}$^{3}$ \\
			
		\vspace{1cm}
		$^1$Department of Industrial Engineering, Bilkent University, Ankara, Turkey \\
		$^2$Department of Industrial Engineering, TED University, Ankara, Turkey \\ 
		$^3$School of Management, University of Bath, Bath, UK
		
		\vspace{24pt}
		\textbf{Correspondence} \\
			Serhat Gul, Department of Industrial Engineering \\ TED University, 06420, \c{C}ankaya, Ankara, Turkey \\ 
			E-mail: serhat.gul@tedu.edu.tr 
		
	\end{center} \vspace{0.2in}
	
\end{titlepage}

\vbox{%
\begin{center}
\large
\textbf{A Stochastic Programming Approach to Surgery Scheduling under Parallel Processing Principle}
\end{center}

\begin{abstract}
	Parallel processing is a principle which enables simultaneous implementation of anesthesia induction and operating room (OR) turnover with the aim of improving OR utilization. In this article, we study the problem of scheduling surgeries for multiple ORs and induction rooms (IR) that function based on the parallel processing principle under uncertainty. We propose a two-stage stochastic mixed-integer programming model considering the uncertainty in induction, surgery and turnover durations. We sequence patients and set appointment times for surgeries in the first stage and assign patients to IRs at the second stage of the model. We show that an optimal myopic policy can be used for IR assignment decisions due to the special structure of the model. We minimize the expected total cost of patient waiting time, OR idle time and IR idle time in the objective function. We enhance the model formulation using bounds on variables and symmetry-breaking constraints. We implement a novel progressive hedging algorithm by proposing a penalty update method and a variable fixing mechanism. Based on real data of a large academic hospital, we compare our solution approach with several scheduling heuristics from the literature. We assess the additional benefits and costs associated with the implementation of parallel processing using near-optimal schedules. We examine how the benefits are inflated by increasing the number of IRs. Finally, we estimate the value of stochastic solution to underline the importance of considering uncertainty in durations. 
	\end{abstract}
	}

{\bf Funding:} This material is based in part on work supported by the TED University Institutional Research Fund [Grant T-18-B2010-33020].


\section{Introduction}
Operating room (OR) expenses account for a significant portion of total costs in a hospital \citep{Gul2015}. \textit{Parallel surgery processing} has recently been implemented in some hospitals to improve OR utilization, and hence reduce costs \citep{Marjamaa2009}. Parallel processing is a technique that allows concurrent implementation of \textit{anesthesia induction} and \textit{turnover}. This implies that the surgical process for a patient can start by giving an anesthetic agent to the patient in the preoperative holding area while the OR is being cleaned and set up for the patient. In other words, it is not mandatory to finish the whole process related to the preceding surgery before giving anesthetic agent to the current patient. Contrary to conventional practice (i.e., \textit{serial processing}), OR is used only for \textit{incision} and closing the incision (i.e., \textit{postincision}) in parallel processing. Therefore, anesthesia induction (i.e., \textit{preincision}) is handled in a room usually called as \textit{induction room} (IR) in the preoperative holding area. 

The benefit of parallel processing has been well acknowledged by the medical community thanks to the actual experiments conducted in several different hospitals. The experiments show that the surgical cases in an OR finish earlier in the  day if parallel processing is preferred over serial processing \citep{Hanss2005, Torkki2005, Harders2006, Friedman2006, Smith2008}. Therefore, the number of surgeries that can be performed in a day increases under this setting. For example,  the  time  without  operation in the ORs  is  decreased  by  45.6\% and 54\%  in the experiments by \citet{Torkki2005} and \citet{Harders2006}, respectively. \citet{El-Boghdadly2020} review 15 studies that involve approximately 8900 patients in the trials in total, and indicate that parallel processing is favored in terms of clinical outcomes. Furthermore, they report that the OR throughput per day increases by 1.7 cases when parallel processing is implemented. \citet{Torkki2005} mention that the benefit of parallel processing would be augmented if surgery schedules are systematically created. However, the optimization models and solution methods that have been developed to solve surgery scheduling problems in the literature have not considered the challenges associated with parallel processing.

Designing surgery schedules is difficult due to uncertainty in surgery durations. In serial processing, surgery duration refers to the summation of preincision, incision and postincision durations \citep{Batun2011}. However, since induction is conducted in a separate location in parallel processing, surgery duration represents the summation of only incision and postincision durations in this study. Uncertainty in induction and OR turnover (i.e., cleaning and set up) times creates additional challenges for the scheduler. Furthermore, the decisions affect the trade-off between patient and provider-related performance measures such as patient waiting time, IR idle time and OR idle time. 

In this article, we study the problem of scheduling surgeries under uncertainty for a surgical suite consisting of multiple ORs and IRs that work based on parallel processing principle. The uncertainty exists due to induction, surgery and turnover durations. We model our problem as a two-stage stochastic mixed integer programming (SMIP) formulation. The main decisions in the model include sequencing patients in the daily appointment list and setting appointment times. We also consider the limited availability of induction rooms. The assignment of patients to induction rooms are the decisions made at the second stage for each realized scenario. We minimize the expected total cost of patient waiting time, OR idle time, and IR idle time in the objective function. We propose a progressive hedging algorithm (PHA) that benefits from the underlying problem structure to design near-optimal schedules. We conduct computational experiments to illustrate the performance of the PHA solutions in comparison to the optimal solutions. The experiments are conducted based on real data of a major hospital in the United States. We compare the proposed method with well-known scheduling heuristics that are used to create daily surgery schedules. We compare parallel processing with serial processing based on near-optimal schedules in terms of patient waiting time and OR closure times. We assess the role of the number of IRs while implementing parallel processing. Finally, we estimate the value of stochastic solution to show the benefit of considering uncertainty in durations while scheduling surgeries under parallel processing principle.  

The organization of the following sections of the article is as follows. In Section 2, the relevant articles are reviewed, and the principal contributions of our study are provided. In Section 3, the problem is described and model formulation is extensively discussed. In Section 4, the proposed PHA algorithm is presented in detail. The results of computational experiments are shown in Section 5. Finally, conclusions and possible extensions are provided in Section 6.

\section{Literature Review}
We review three categories of articles to distinguish our study from the other studies in the literature: studies (i) considering appointment scheduling decisions for a daily list of patients, (ii) focusing on benefits obtained through concurrent implementation of anesthesia induction and turnover by parallel processing principle, and (iii) defining parallel processing as the setting where a single surgeon works in two different ORs in an alternating fashion.

Among the first category of studies, \citet{Gupta2008} conduct a literature review on appointment time scheduling problems and group the studies based on the settings considered. These settings are primary care, specialty care, and elective surgery care. We only present studies that consider elective surgery care. Studies of the remaining two groups are extensively discussed in \citet{Cayirli2003}, \citet{Gupta2008} and \citet{AhmadiJavid2017}. 
\citet{Gul2011} use heuristics to sequence surgeries and set patient appointment times with the objective of minimizing expected patient waiting time and surgical suite overtime. For setting patient appointment times, they use a job hedging procedure, where the allocated duration for a surgery is calculated based on a certain percentile of the surgery duration distribution. 

\citet{Denton2003} formulate the appointment scheduling problem for a single-server case as a two-stage stochastic linear program with the objective of minimizing the expected total cost of customer waiting time, server idle time, and tardiness. Their model is generalized to represent any service system where appointment scheduling would be relevant. However, their study is motivated by a surgery scheduling problem in a single OR. In the context of surgery scheduling, the performance measures correspond to patient waiting time, OR idle time, and OR overtime, respectively. They solve the resulting formulation using an L-shaped algorithm with sequential bounding. \citet{Khaniyev2020} focus on the same problem and propose simple heuristics that are shown to perform well for a daily sequence of surgeries with non-identical duration distributions. 

\citet{Denton2007} extend the model in \citet{Denton2003} by also considering sequencing decisions. They formulate the problem as a two-stage stochastic mixed integer program. They make sequencing decisions based on heuristics and solve the reduced version of the problem to set patient appointment times. \citet{Mancilla2012} provide a reformulation of the model in \citet{Denton2007}, and implement a heuristic based on Benders’ decomposition to solve the model.

\citet{Zhang2015} study how to set appointment times for surgeries in a multiple-OR setting. The objective is to minimize the total expected cost of surgeon waiting, OR idle, and OR overtime. They formulate a simulation-based optimization model. They propose stochastic gradient algorithm to solve the model.\citet{Vandenberghe2019} study a stochastic surgery sequencing problem also for a multiple-OR case by assuming that the OR assignment decisions are already made. They call their problem as the stochastic break-in moment problem since they minimize the expected maximum waiting time of an arriving emergent surgery. They solve the model using sample average approximation algorithm and different heuristics. \citet{Atighehchian2020} first allocate surgeries to ORs, surgery assistants and residents by considering the setting of teaching hospitals. They then sequence surgeries for each OR and resident without setting any appointment times. The objective of the model minimizes expected OR idle time and overtime costs. They solve the resulting two-stage stochastic programming model using an L-shaped method.  

\citet{Lee2014} set patient appointment times for multiple ORs by considering the availability of post-anesthesia care unit (PACU) resources and assuming a fixed patient sequence. The problem is formulated as flexible job shop model with fuzzy sets which are used to model the uncertainty in process durations. The criteria include patient waiting time, OR idle time, and surgical suite completion time in the study. They use a combination of heuristic and genetic algorithms to find solutions. \citet{Bai2017} also model the appointment scheduling problem for multiple ORs considering PACU constraints. Their model also does not allow changes in the order of surgeries. They formulate the problem as a discrete-event dynamic system-based stochastic optimization model. They minimize the expected total cost of patient waiting time, OR idle time, OR blocking time, OR overtime and PACU overtime. To solve the deterministic equivalent of their model, they propose a gradient based algorithm.   \citet{Varmazyar2020} study patient sequencing problem for a setting also including multiple ORs and PACU resources. They estimate surgery and PACU durations using continuous phase-time distributions and propose a heuristic to solve the problem. 

\citet{Erdogan2013} formulate a dynamic appointment scheduling problem using a multi-stage stochastic linear program. They schedule appointment requests dynamically. Their study is motivated by the arrival of urgent add-on cases that require revising the current schedule dynamically. They solve the model using decomposition-based algorithms. 

Our study focuses on appointment scheduling decisions as the first category of articles reviewed. However, we consider the flow of patients through a surgical suite including multiple IRs and ORs that operate based on parallel processing principle.   

In the second category of our literature review, we only identify the study of \citet{Marjamaa2009}. They investigate four different parallel processing systems in which the induction and turnover are carried out simultaneously using operations research methods. They compare their models with serial processing models using discrete-event simulation.  All models are assumed to have four ORs. In the first model, there is an IR adjacent to each OR. In this model, there is an additional induction team consisting of two nurses for each OR. While a team is operating in the OR, the other team with an anesthesiologist deals with the patient in the IR.  In the second model, the team consisting of an anesthesiologist and two nurses induce anesthesia to all patients in different IRs. In the third model, there is a central IR, where there are three beds and an induction team. In the last model, there are three ORs used for surgeries and one OR used for induction. The results are provided in terms of the total number of completed surgeries, over-utilization and under-utilization of ORs, and personnel expenses. All parallel processing systems yield lower costs than serial processing systems. The study also investigates the possible benefit of assigning all short surgeries to one OR and the long surgeries to another OR. 

\citet{Marjamaa2009} investigate parallel processing systems only using simulation. On the other hand, we develop a two-stage SMIP model to find near-optimal schedules for surgical suites that operate based on parallel processing.

The third category of studies investigate a different form of parallel processing where a single surgeon performs parallel surgeries in two ORs.  Once the incision process of the surgery is completed, the surgeon moves to the other OR and starts dealing with the incision phase of that patient. However, the anesthesia inductions are applied in the ORs. \citet{Batun2011} and \citet{Mancilla2013} are the two studies that propose such parallel processing implementations. 

\citet{Batun2011} decide on the number of ORs that would be used during the day, assign surgeries to ORs, sequence surgeries, and set the start time of the first surgery of each surgeon in a day.  It is assumed that the surgeons operate their surgeries in series. Hence, the appointment time for a particular surgery is not needed to be set. The problem is modeled as an SMIP and solved by an L-shaped algorithm. \citet{Mancilla2013} study the assignment of surgeries to ORs and sequence surgeries. The problem is modeled as an SMIP and solved by a decomposition based method.

In our study, parallel processing represents a setting which is completely different from the one considered in \citet{Batun2011} and \citet{Mancilla2013}. The differences of our study can be summarized as follows: An induction cannot be applied in an OR. An IR, located adjacent to an OR in the preoperative holding area, is used for this task instead. The operating discipline is considered to be parallel because the induction for a patient in an IR and the turnover after the surgery of the preceding patient in an OR are performed simultaneously. A surgeon works only in a single OR and does not have to alternate between ORs for the implementation of parallel processing. Furthermore, surgeons do not have to operate their own surgeries consecutively. Therefore, we set appointment times for each surgery in our study.

\section{Problem Description and Model Formulation}

We study the problem of sequencing patients and setting appointment times for a daily surgery list of a surgical suite having multiple ORs and IRs that function according to the parallel processing principle. Therefore, the induction is not applied in an OR, but in an IR so that the patient can be ready for the incision part of the surgery while the OR is cleaned and set up after the previous surgery. We assume that patients are punctual and arrive to the surgical suite exactly at their appointment times. They first wait for an IR to become available if all are busy. Then, they can be assigned to any one of the available IRs, as the IRs are assumed to be identical. After the induction is applied by an anesthesiologist in the IR, the patient is taken to an OR for surgery, possibly after some waiting. We ignore the transfer time from an IR to an OR. No OR assignment decision is made in the problem since we assume that patient-to-OR assignments are known a priori. After the surgery is performed in the OR, the turnover starts to make the OR ready for the next patient.  

We assume that an anesthesiologist provides care only for a set of patients whose surgeries are planned to be performed in the same OR. Since the OR assignments are already known, this means the anesthesiologist of each patient is also known. Furthermore, we assume the anesthesiologist can start the induction for a patient in an IR even if the surgery of the preceding patient continues in the OR that the anesthesiologist is responsible for. This is because the anesthesiologist does not need to be active during the whole surgery process of a patient as their aides can also monitor the patient in the OR. On the other hand, we assume that an anesthesiologist must be present in the IR over the whole induction process for a patient.

We formulate the problem as a two-stage SMIP. At the first stage of the model, patients are sequenced, and their appointment times are set under uncertainty in induction, surgery and turnover durations. At the second stage, patients are assigned to one of the identical IRs based on a \textit{myopic policy} that is optimal due to the special structure of the two-stage SMIP (as will be explained after the model is provided). The myopic policy allows all IR assignment decisions to be made at a single stage. The performance measures including IR and OR idle times, and patient waiting times for IR and OR are also calculated in the second stage. Besides, some auxiliary variables including surgery start times, IR closure times and OR closure times are also calculated at the same stage. 

The waiting time of a patient for an IR is equal to the difference between the induction start time and appointment time of the patient. The induction start of a patient may be delayed because the induction of preceding patient may not finish on time. In that case, the induction of the patient can start only when the preceding patient is taken to the OR for surgery. This implies that the patients stay in an IR after induction when they need to wait for an OR to become available. The waiting time of a patient for an OR is equal to the difference between the surgery start time and induction finish time of the patient. The surgery of a patient may start late because the turnover for the preceding surgery may finish later than planned.   

Idle time of an IR is equal to the difference between the IR closure time and summation of induction durations of all patients treated in that IR. Note that an IR is assumed to be idle even if a patient occupies the room as he/she waits for an OR to be available, because the IR staff is kept idle during that period. This implies that such idleness cases are punished twice due to their negative impact into both provider and patient related performance measures. An IR closes when the last patient treated in the IR is taken from the IR to an OR for surgery. Therefore, IR closure time is equal to the surgery start time of the last patient treated in the IR. Idle time of an OR is equal to the difference between the OR closure time and summation of surgery and turnover durations of all patients treated in the OR. An OR closes when the turnover for the last patient treated in the OR is finished. 

\begin{figure}
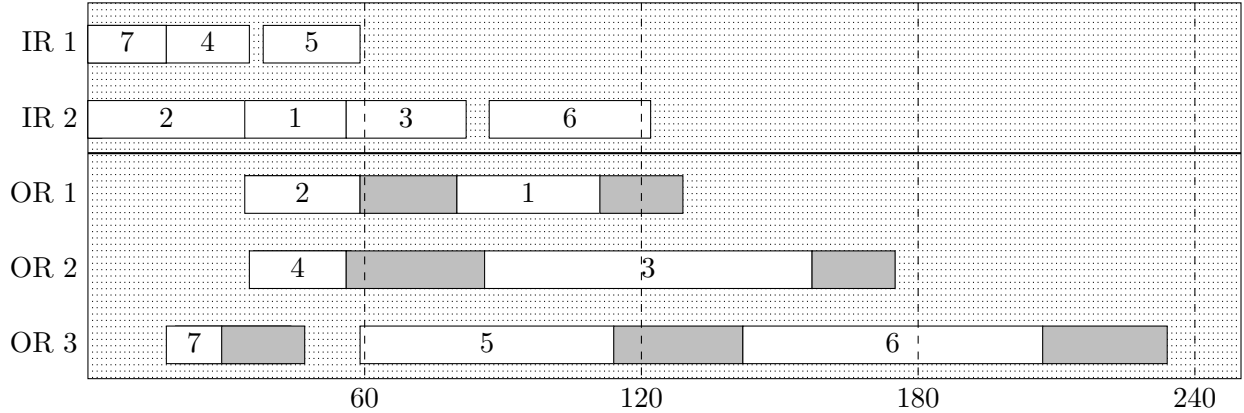

\begin{ganttchart}[vgrid,,expand chart=\textwidth, bar height = 0.5]{1}{250}
\ganttbar{IR 1}{1}{17}
\ganttbar[inline]{7}{1}{17}
\ganttbar[inline]{4}{18}{35}
\ganttbar[inline]{5}{39}{59}\\
\ganttbar{IR 2}{1}{3} 
\ganttbar[inline]{2}{1}{34}
\ganttbar[inline]{1}{35}{56}
\ganttbar[inline]{3}{57}{82}
\ganttbar[inline]{6}{88}{122}
\ganttnewline[thick]

\ganttbar{OR 1}{35}{59}
\ganttbar[inline]{2}{35}{59}
\ganttbar[inline, bar/.append style={fill=gray!50}]{}{60}{80}
\ganttbar[inline]{1}{81}{111}
\ganttbar[inline, bar/.append style={fill=gray!50}]{}{112}{129}\\
\ganttbar{OR 2}{37}{57}
\ganttbar[inline]{4}{36}{56}
\ganttbar[inline, bar/.append style={fill=gray!50}]{}{57}{86}
\ganttbar[inline]{3}{87}{157}
\ganttbar[inline, bar/.append style={fill=gray!50}]{}{158}{175}\\
\ganttbar{OR 3}{20}{44}
\ganttbar[inline]{7}{18}{29}
\ganttbar[inline, bar/.append style={fill=gray!50}]{}{30}{47}
\ganttbar[inline]{5}{60}{114}
\ganttbar[inline, bar/.append style={fill=gray!50}]{}{115}{142}
\ganttbar[inline]{6}{143}{207}
\ganttbar[inline, bar/.append style={fill=gray!50}]{}{208}{234}

\ganttvrule[vrule/.append style={thin}]{60}{60}
\ganttvrule[vrule/.append style={thin}]{120}{120}
\ganttvrule[vrule/.append style={thin}]{180}{180}
\ganttvrule[vrule/.append style={thin}]{240}{240}
\end{ganttchart}
\caption{An example 4-hour IR and OR schedule, where each number corresponds to a patient, white bars to the duration of induction and surgeries, and dark bars to turnover durations}
\label{fig:schedule}
\end{figure}

To illustrate the flow of operations in the surgical suite and explain the performance measures clearly, Figure \ref{fig:schedule} shows an example 4-hour schedule with seven patients, two IRs, and three ORs. Patients 1 and 2 are pre-assigned to OR 1, patients 3 and 4 to OR 2, and the remaining patients to OR 3. The appointment times are $a_1=20$, $a_2=a_7=0$, $a_3=55$, $a_4=10$, $a_5=38$, and $a_6=87$. In the beginning of the shift, patients 2 and 7 are assigned to IRs 2 and 1, respectively. Patient 4 arrives during the induction of these two patients, and has to wait for 7 minutes until patient 7 finishes induction and moves to OR 3. Similarly, patient 1 arrives while the induction of patients 2 and 4 are in progress, and needs to wait for 14 minutes before the induction of patient 2 ends and surgery starts in OR 1. One minute later, the induction of patient 4 finishes and the patient moves to OR 2. Patients 5, 3 and 6 arrive later in order, and the only other waiting for the IR is patient 3, for 1 minute.

Figure \ref{fig:schedule} also shows examples of parallel processing in the suite. For example, the induction of patient 5 is underway while the turnover in its assigned OR is in progress after the surgery of patient 7 is complete. A similar observation can be made for patient 3 in IR 2, at the same time the turnover in OR 2 is ongoing prior to the surgery of this patient.

For the example in Figure \ref{fig:schedule}, the total waiting time for the IR is 22 minutes (14, 1, and 7 minutes for patients 1, 3, and 4, respectively). Patients need to wait for a total of 48 minutes for the ORs (24, 4, and 20 minutes for patients 1, 3, and 6, respectively). This yields a total waiting time of 70 minutes. IR 1 closes after 59 minutes, when patient 1 moves to OR 1. IR 2 closes after 142 minutes, after OR 3 becomes available for patient 6. Similarly, ORs 1, 2 and 3 close after 142, 175 and 234 minutes, respectively, following the turnover of their last patients. Total OR idle time is 98 minutes, with 34, 35 and 29 minutes of idle time at OR 1, OR 2, and OR 3, respectively. IR 1 has a total idle time of 4 minutes between the departure of patient 4 for OR 2 and the arrival of patient 5. IR 2 has a total idle time of 25 minutes; 5 minutes between the departure of patient 3 for OR 2 and the arrival of patient 6, and 20 minutes when patient 6 needs to wait for the surgery of patient 5 to end in OR 3.

\begin{table} 
\centering
\caption{Notation used throughout the mathematical model}
\label{tab:notation}
\vspace{10pt}
\label{tab:notation}
\begin{tabular}{cl}
\multicolumn{2}{l}{\textbf{Index sets}} \\
\hline
$S$ & Surgeries (patients) \\
$K$ & Induction rooms \\
$R$ & Operating rooms \\
$S_r$ & Surgeries to be conducted in OR $r \in R$ \\
$\Omega$ & Scenarios \\
\hline
& \\
\multicolumn{2}{l}{\textbf{Deterministic Parameters}} \\
\hline
$c^I$ & Cost of idle OR time per minute \\
$c^J$ & Cost of idle IR time per minute \\
$c^W$ & Cost of patient waiting time per minute \\
$M$ & A large value\\
\hline
& \\
\multicolumn{2}{l}{\textbf{Stochastic Parameters}} \\
\hline
$d_i(\omega)$ & Duration of surgery for patient $i\in S$ under scenario $\omega \in \Omega$ \\
$e_i(\omega)$ & Duration of induction for patient $i\in S$ under scenario $\omega \in \Omega$ \\
$q_i(\omega)$ & Duration of turnover after the surgery for patient $i\in S$ under scenario $\omega \in \Omega$ \\
\hline
& \\
\multicolumn{2}{l}{\textbf{First-Stage Decision Variables}} \\
\hline
$a_i$ & Appointment time of patient $i\in S$ \\
$u_{ij}$ & $=\begin{cases} 1, & \text{if patient } i \in S \text{ precedes patient } j\in S \text{ (general precedence)} \\ 0, & \text{otherwise} \end{cases}$ \\\hline
& \\
\multicolumn{2}{l}{\textbf{Second-Stage Decision Variables}} \\
\hline
$y_{ik}(\omega)$ & $=\begin{cases} 1, & \text{if patient } i \in S \text{ is assigned to IR } k\in K \text{ under scenario } \omega \in \Omega \\ 0, & \text{otherwise} \end{cases}$ \\
$A_i(\omega)$ & Surgery start time for patient $i\in S$ under scenario $\omega \in \Omega$ \\
$W_i(\omega)$ & Waiting time for patient $i\in S$ for an OR under scenario $\omega \in \Omega$ \\
$Y_i(\omega)$ & Waiting time for patient $i\in S$ for an IR under scenario $\omega \in \Omega$ \\
$I_r(\omega)$ & Idle time of OR $r\in R$ under scenario $\omega \in \Omega$ \\
$J(\omega)$ & Total idle time of all IRs under scenario $\omega \in \Omega$ \\
$F_r(\omega)$ & Closure time of OR $r\in R$ under scenario $\omega \in \Omega$ \\
$G_k(\omega)$ & Closure time of IR $k\in K$ under scenario $\omega \in \Omega$ \\
\hline
\end{tabular}
\end{table}

The notation we use for the formulation of the problem is provided in Table \ref{tab:notation}. Using the sets, parameters and decision variables in the table, we formulate the problem using the following two-stage SMIP:
\begin{alignat}{2}
\min \quad & \mathcal{Q}(\textbf{a,u}) \quad && \label{1a}\\
s.t. \quad  & a_j  \geq a_i - M (1- u_{ij}) \quad && \forall{i,j} \in S, j\neq{i} \label{1b}\\
& u_{ij}+u_{ji}=1 \quad && \forall{i,j} \in S, j>i  \label{1c} \\
& a_i \in \mathbb{Z}^+ && \forall{i} \in S  \label{1d} \\
& u_{ij}\in {\{0,1\}} \quad && \forall{i,j} \in S, j\neq{i} \label{1e}
\end{alignat}
where $\mathcal{Q}(\textbf{a,u})=E_{\xi}[Q(\textbf{a,u},\xi(\omega))$ is the expected second-stage cost, and $Q(\textbf{a,u},\xi(\omega)))$ is determined by:
\begin{alignat}{2}
\min \quad & c^JJ(\omega)+c^I\sum_{r\in{R}}I_r(\omega)+c^W\sum_{i\in{S}}\big(W_i(\omega)+Y_i(\omega) \big ) \quad && \label{2a}\\
s.t. \quad  & \sum_{k\in{K}}y_{ik}(\omega)=1 \quad && \forall i \in S  \label{2d}\\
& a_j+Y_j(\omega)+M(1-u_{ij})\ge a_i+Y_i(\omega) \quad && \forall{i,j} \in S, j\neq{i} \label{2e} \\
& a_j+Y_j(\omega) +M(1-u_{ij})+M(2-y_{ik}(\omega)-y_{jk}(\omega))\ge{A_i(\omega)} \quad && \forall{i,j} \in S, j\neq{i}, \forall k \in K  \label{2f} \\
& a_j+Y_j(\omega) + M(1- u_{ij}) \geq a_i+Y_i(\omega) + e_i(\omega) \quad && \forall{i,j} \in{S_r},j\neq{i}, \forall r \in R \label{2g} \\
& A_i(\omega)= a_i+Y_i(\omega)+ e_i(\omega) +W_i(\omega) \quad && \forall i \in S \label{2h} \\
& A_j(\omega)+M(1-u_{ij})\ge A_i(\omega)+d_i(\omega) +q_{i}(\omega) \quad && \forall{i,j} \in{S_r},j\neq{i}, \forall r \in R \label{2i} \\
& F_r(\omega)\ge A_i(\omega)+d_i(\omega)+q_{i}(\omega) \quad && \forall i\in{S_r}, \forall r \in R  \label{2k} \\
& I_r(\omega)={F_r(\omega)-\sum_{i\in{S_r}}(d_i(\omega)+q_i(\omega))} \quad && \forall {r} \in R  \label{2l} \\
& G_k(\omega)\ge{A_i(\omega)-M(1-y_{ik}(\omega)}) \quad && \forall i \in S, \forall k \in K \label{2m} \\
& J(\omega)=\sum_{k\in{K}}G_k(\omega)-\sum_{i\in{S}}e_i(\omega) \quad && \label{2n} \\
& A_i(\omega),W_i(\omega),Y_i(\omega)\ge{0} \quad && \forall i \in S \label{2o} \\
& I_r(\omega),F_r(\omega),G_k(\omega),J(\omega)\ge{0} \quad && \forall r \in R, \forall k \in K  \label{2p} \\
& y_{ik}(\omega) \in{\{0,1\}} \quad && \forall i \in S, \forall k \in K \label{2r}
\end{alignat}

\noindent Assuming that the uncertainty in durations can be represented by a finite set of scenarios, the objective function \eqref{1a} only minimizes the expected second-stage cost over all scenarios, which is equivalent to the total expected cost of patient waiting for IR and OR, IR idle time, and OR idle time. 
 
The first-stage constraints are shown by the constraint set \eqref{1b}-\eqref{1e}. Constraints \eqref{1b} and \eqref{1c} define the sequencing rules for patients. Constraint \eqref{1b} ensures that if patient $i \in S$ precedes patient $j \in S$ in the daily patient list, then the appointment time of patient $i$ cannot be set later than that of patient $j$. Constraint \eqref{1c} provides that either patient $i\in S$ precedes patient $j \in S$ in the daily patient list or vice versa. Constraints \eqref{1d} and \eqref{1e} represent integrality and binary restrictions on the first-stage variables, respectively.

The second-stage objective function for given first-stage decision variables and a particular scenario is formulated explicitly by \eqref{2a}. It minimizes the total cost of patient waiting for IR and OR, IR idle time, and OR idle time for a given scenario. 

Constraint \eqref{2d} ensures that each patient is assigned to exactly one IR. Constraint \eqref{2e} ensures that the induction start time of patient $j \in S$ is greater than or equal to that of patient $i \in S$ only in case patient $j$ comes after patient $i$ in the patient sequence. Constraint \eqref{2f} applies only to the patients that are assigned to the same IR. The constraint guarantees that the induction for patient $j \in S$ can start only after patient $i \in S$ is moved from the IR and placed to her OR in case patient $j$ comes after patient $i$ in the patient sequence. 

Constraint \eqref{2g} enforces a restriction for the patients who are pre-assigned to the same OR. Due to one-to-one matching between ORs and anesthesiologists, the constraint is valid for the patients who are pre-assigned to the same anesthesiologist. For a given set of surgeries to be treated in OR $r \in R$, it ensures that the induction for patient $j \in S_r$ can start after that for patient $i \in S_r$ ends, if patient $i$ precedes patient $j$ in the patient sequence. The constraint implies that each anesthesiologist can perform a single induction at a time even if the patients are treated in separate IRs.

Constraint \eqref{2h} is formulated to determine the amount of patient waiting time for an OR right after the induction ends. The constraint sets the surgery start time for patient $i \in S$ equal to the summation of induction completion time of patient $i$ and waiting time of patient $i$ for an OR. Constraint \eqref{2i} ensures that surgery for a patient can start only if the turnover of the earlier surgeries in the same OR ends. The constraint also implies that the patient sequence in the first stage is preserved at the second stage for each OR. Constraint \eqref{2k} sets the actual closure time of an OR. It ensures that the closure time of OR $r \in R$ is not earlier than the turnover finish time of any surgery performed in OR $r$. Constraint \eqref{2l} calculates the idle time in an OR. OR idle time is set equal to the actual OR closure time minus the total surgery and turnover durations for the patients treated in the OR. Constraint \eqref{2m} defines the closure time of an IR. The closure time of an IR cannot be earlier than the surgery start times of the patients who go through that IR, as the patients wait for their surgeries in the IRs. Constraint \eqref{2n} calculates the total idle time over all IRs. This is achieved by subtracting the summation of induction durations from the summation of IR closure times. Finally, constraints \eqref{2o}-\eqref{2r} define non-negativity and binary restrictions on the second-stage variables. 

We consider IR assignment decisions in the second stage, because considering them in the first stage would result in inefficiencies in the surgical suite as the IRs are identical. Even though the IR assignment decisions are made dynamically in the surgical suites as patients arrive for their appointments, the special structure of our model allows us to avoid formulating the problem as a multi-stage SMIP. We show that after the sequence of patients is fixed in the first stage, the decisions in the second stage can be made based on a myopic policy. Since using a myopic policy prevents the need of considering uncertainty related to induction, surgery and turnover durations, our problem can be formulated as a two-stage SMIP. The justification of our modeling choice is based on the following proposition, whose proof is provided in Appendix A.

\begin{proposition}
	For a given sequence of patients obtained from the first stage, each arriving patient is assigned to the first available induction room in the optimal schedule.
\end{proposition}

The proof of proposition is valid mainly due to constraints (\ref{2e}). These constraints ensure that the sequence of patients in the first stage does not change during the treatments of the patients both in the IR and OR phases in the second stage. Hence, the optimal assignment of the patients to IRs can be achieved by a non-idling policy. This implies that all assignment decisions can be made at a single stage without formulating a dynamic model. 

Note that the assumption of unchanged patient sequence is also practically relevant. Since we consider only elective surgeries in our problem, the urgency levels of the surgeries are all the same. Furthermore, the ORs to be used for each surgery are predetermined. Therefore, switching the planned position of two patients waiting for their surgeries would not be realistic considering the possible reactions of the patients. Such a change would be perceived by a patient as an unfair treatment, and hence the satisfaction level of the patients would be significantly affected in a negative direction. Therefore, the managers instead must ensure that the appointment times are carefully set such that the smooth flow of patients through the surgical suite can be maintained without making a change in the original sequence. 

We strengthen the two-stage SMIP formulation by adding a number of valid inequalities into the second-stage problem, and solve the resulting formulation in the computational experiments. In particular, we impose bounds on some variables and formulate symmetry-breaking constraints.

We first propose a bound on the variables representing waiting times of patients for an OR. Suppose that patients $i, j, k\in S_r$  are the first, second and third patients to be treated in OR $r \in R$. The worst-case scenario from the perspective of patient $j$ and $k$ occurs when the induction finish times of patient $j$ and $k$ are equal to the surgery start time of patient $i$. Since patient $j$ must wait until the OR is made ready after patient $i$ leaves, the maximum amount that patient $j$ waits for an OR is equal to the summation of surgery and turnover durations for patient $i$. Moreover, the maximum amount that patient $k$ waits for an OR is equal to the summation of surgery and turnover durations for patients $i$ and $j$. By generalizing this relationship for the patients to be treated in the same OR, we formulate the following bound on the $W_j(\omega)$ variables:   
\begin{flalign}
&W_j(\omega)\le{\sum_{i\in{S_r},i\ne{j}}}\big(d_i(\omega)+q_i(\omega)\big)u_{ij}\qquad  \forall{j}\in{S_r},\forall r \in R \label{ub}
\end{flalign}

Next, we enforce that the closure time of the IR with index  $(k-1) \in K$ is greater than or equal to the closure time of the IR with index $k \in K, k>1$ using the following constraint: 
\begin{flalign}
&G_{k-1}(\omega)\ge{G_{k}(\omega)}\qquad \forall k \in K \setminus \{1\} \label{sym2}
\end{flalign}
\section{Solution Methodology}

The number of variables and constraints would be quite large for realistically sized instances of SMIP models, since there are as many second-stage subproblems as the number of scenario realizations. The resulting complexity necessitates the application of sophisticated algorithms even when a near-optimal solution of a model is sufficient \citep{birge2011introduction}. Therefore, various decomposition algorithms have been proposed to solve SMIP models. The L-shaped algorithm of \citet{VanSlykeWets1969} is not applicable to the models including integer variables, which is also the case in our study. \citet{laporte1993integer} propose a variant of this algorithm, called integer L-shaped algorithm, to solve SMIP models having only binary variables in the first stage. Value function and set convexification based approaches are also developed in the literature to deal with models having binary variables in both stages \citep{sen2005c3,yuan2009enhanced}. However, they cannot be utilized to solve our model, because our model includes both general and binary integer variables in the first stage. Having also mixed-binary variables in the second stage makes the solution process particularly difficult for our model.  

A class of algorithms that would be appropriate for solving our model is scenario decomposition-based methods. The PHA, which was proposed by \citet{Rockafellar1991}, is a commonly used scenario decomposition algorithm \citep{Hvattum2009, Goncalves2012,Gul2015,demiretal21}. The algorithm is basically a type of augmented Lagrangian relaxation approach. To apply scenario decomposition at each iteration, a reformulation of the model is required before the implementation of the algorithm. Next, the nonanticipativity constraints (i.e., the constraints that enforce the restrictions that the first-stage variable values cannot take different values under different scenarios) in the reformulated model are relaxed and carried to the objective function using a penalty and Lagrangian term. The augmented Lagrangian relaxation approach allows decomposing the model into scenario subproblems. After solving scenario subproblems independently and then aggregating the scenario solutions over multiple iterations, the PHA provides a nonanticipative and feasible solution for the original model.   

The PHA yields optimal solutions for convex models \citep{Rockafellar1991}. Since the two-stage SMIP model is nonconvex due to binary variables at both stages of the formulation, the algorithm provides a heuristic solution to our problem. Note that the PHA is suggested as a good heuristic for SMIP models in several studies in the literature \citep{Crainic2011, Watson2011,Gul2015, demiretal21}.

We first provide the reformulation of the two-stage SMIP model in the remaining parts of this section. Next, we discuss the basic PHA with its general steps. Finally, we present our own implementation of the PHA which improves the performance of the basic PHA by linearizing the objective function, and using a variable fixing method and a penalty update mechanism.
\subsection{Problem reformulation} 
To implement the PHA, a scenario separable formulation is needed. Therefore, we first reformulate the two-stage SMIP by explicitly representing nonanticipativity constraints in the model. For this end, we additionally define two decision variables: $a_i(\omega)$ as the appointment time for  patient $i \in S$ under scenario $\omega \in \Omega$ and $u_{ij}(\omega)$ as a binary variable denoting whether patient $i \in S$ precedes patient $j \in S$ or not under scenario $\omega \in \Omega$.

Denoting the probability of occurrence of each scenario $\omega \in \Omega$ by $p^\omega$, the reformulation of the two-stage SMIP, which we refer to as R-SMIP, is obtained as follows:
{\footnotesize
\begin{alignat}{2}
\min \quad & \sum_{\omega\in \Omega} p^\omega\Big(c^JJ(\omega)+c^I\sum_{r\in R}I_r(\omega)+c^W\sum_{i\in S}\big(W_i(\omega)+Y_i(\omega)\big)\Big) \quad && \label{ra}\\
s.t. \quad  & a_j(\omega)  \geq a_i(\omega) - M (1- u_{ij}(\omega)) \quad && \forall{i,j} \in S, j\neq{i}, \forall \omega \in \Omega \label{r1}\\
& u_{ij}(\omega)+u_{ji}(\omega)=1 \quad && \forall{i,j} \in S, j>i, \forall \omega \in \Omega \label{r2} \\
& \sum_{k\in{K}}y_{ik}(\omega)=1 && \forall i \in S, \forall \omega \in \Omega  \label{r5} \\
& a_j(\omega)+Y_j(\omega)+M(1-u_{ij})\ge a_i(\omega)+Y_i(\omega) \quad && \forall{i,j} \in S, j\neq{i},\forall \omega \in \Omega \label{r6} \\
& a_j(\omega)+Y_j(\omega)+M(1-u_{ij})+M(2-y_{ik}(\omega)-y_{jk}(\omega))\ge{A_i(\omega)} \quad && \forall{i,j} \in S, j\neq{i},\forall k \in K, \forall \omega \in \Omega \label{r7} \\
& a_j(\omega)+Y_j(\omega) + M(1- u_{ij}) \geq a_i(\omega)+Y_i(\omega) + e_i(\omega) \quad && \forall{i,j}\in{S_r},j\neq{i}, \forall r \in R, \forall \omega \in \Omega\label{r8} \\
& A_i(\omega)=a_i(\omega)+Y_i(\omega)+e_i(\omega)+W_i(\omega) \quad && \forall i \in S, \forall \omega \in \Omega \label{r9} \\
& A_j(\omega)+M(1-u_{ij})\ge A_i(\omega)+d_i(\omega)+q_{i}(\omega) \quad && \forall {i,j} \in{S_r},j\neq{i}, \forall r \in R, \forall \omega \in \Omega \label{r10} \\
& F_r(\omega)\ge A_i(\omega)+d_i(\omega)+q_{i}(\omega) \quad && \forall i\in{S_r}, \forall r \in R, \forall \omega \in \Omega  \label{r12} \\
& I_r(\omega)={F_r(\omega)-\sum_{i\in{S_r}}(d_i(\omega)+q_i(\omega))} \quad && \forall r \in R, \forall \omega \in \Omega \label{r13} \\
& G_k(\omega)\ge{A_i(\omega)-M(1-y_{ik}(\omega)}) \quad && \forall i \in S, \forall k \in K, \forall \omega \in \Omega \label{r14} \\
& J(\omega)=\sum_{k\in{K}}G_k(\omega)-\sum_{i\in{S}}e_i(\omega) \quad && \forall \omega \in \Omega \label{r15} \\
& a_i(\omega) = a_i \quad && \forall i \in S,\ \forall \omega \in \Omega \label{r16} \\
& a_i: integer \quad && \forall i \in S  \label{r17} \\
& a_i(\omega): integer \quad && \forall i \in S , \forall \omega \in \Omega \label{r18} \\
& u_{ij}(\omega)\in {\{0,1\}} \quad && \forall i,j \in S, j\neq i,\forall \omega \in \Omega  \label{r19} \\
& A_i(\omega),W_i(\omega),Y_i(\omega)\ge{0} \quad && \forall i \in S, \forall \omega \in \Omega \label{r20} \\
& I_r(\omega),F_r(\omega),G_k(\omega),J(\omega)\ge{0} \quad && \forall r \in R, \forall k \in K, \forall \omega \in \Omega  \label{r21} \\
& y_{ik}(\omega) \in{\{0,1\}} \quad && \forall i \in S, \forall k \in K, \forall \omega \in \Omega \label{r22}
\end{alignat}
}%
The main difference of R-SMIP from the two-stage SMIP is that the first-stage variables are defined for each scenario. However, allowing those variables to take different values in different scenarios implies that the scenario realizations in the future would be totally anticipated. Therefore, constraint (\ref{r16}), which maintains nonanticipativity by enforcing that $a_i(\omega)$ take the same values for each $\omega \in \Omega$, is formulated. Note that $a_i$ variables are called as \textit{consensus variables} in R-SMIP. Since  $a_i(\omega)$ values also determine $u_{ij}(\omega)$ values through constraint (\ref{r1}), there is no need to add nonanticipativity constraints explicitly for $u_{ij}(\omega)$ variables.  

The next step towards reaching a separable formulation requires relaxation of the nonanticipativity constraints by penalizing the violation of them in the objective function. The  penalization is achieved through Lagrangian multipliers, $\mu_i(\omega), \forall i \in S, \omega \in \Omega$, and penalty parameter, $\rho$. Note that ordinary Euclidean norm is considered in the term where $\rho$ is used. The formulation obtained after applying augmented Lagrangian relaxation, which we call as SQP, is as follows: 
{\allowdisplaybreaks \begin{flalign}
		&\mbox{min}\sum_{\omega\in \Omega} p^\omega \Big(c^JJ(\omega)+c^I\sum_{r\in R}I_r(\omega)+c^W\sum_{i\in S}\big(W_i(\omega)+Y_i(\omega)\big) \notag \\ 
		& + \sum_{i\in S} \mu_i(\omega) \big(a_i(\omega) - a_i\big) + \frac{\rho }{2} \sum_{i\in S} \| a_i(\omega) - a_i\|^2  \Big)\label{Laga}\\
		&\mbox{s.t.} \nonumber \\
		& (\ref{r1})- (\ref{r15})\\
		& (\ref{r17})- (\ref{r22})
\end{flalign}}
The terms with Lagrangian multipliers and penalty components prevent separability of the formulation due to having consensus variables in them. When consensus variables are replaced by their estimated values, the obstacle can be overcome. For this purpose, a proximal point method is applied by taking weighted sum of consensus variable values for each patient over all scenarios.  The resulting parameter, which we call as \textit{consensus parameter}, is calculated as shown below:  
\begin{equation}
\hat{a}_i = \sum_{\omega\in \Omega} p^\omega a_i(\omega), \quad \forall i \in S
\end{equation}    

The consensus parameter must be updated at each iteration of the PHA. If the consensus parameter value is feasible for R-SMIP, then it would also represent feasible $a_i$ values for the two-stage SMIP. The PHA aims to improve these feasible consensus parameter values over iterations and yield a near-optimal solution for the two-stage SMIP.

Replacing consensus variables in SQP with consensus parameters, we obtain the following separable stochastic quadratic formulation named SSQP.
{\allowdisplaybreaks \begin{flalign}
		&\mbox{min}\sum_{\omega\in \Omega} p^\omega \Big(c^JJ(\omega)+c^I\sum_{r\in R}I_r(\omega)+c^W\sum_{i\in S}\big(W_i(\omega)+Y_i(\omega)\big) \notag \\ 
		& + \sum_{i\in S} \mu_i(\omega) \big(a_i(\omega) - \hat{a}_i\big) + \frac{\rho }{2} \sum_{i\in S} \| a_i(\omega) - \hat{a}_i\|^2  \Big)\label{sa}\\
		&\mbox{s.t.} \nonumber \\
		& (\ref{r1})- (\ref{r15})\\
		& (\ref{r18})- (\ref{r22})
\end{flalign}}
Finally, by applying scenario decomposition on SSQP, we formulate the following scenario subproblems for each scenario $\omega \in \Omega$ , called SSP, to be solved independently at each iteration of the PHA:
\begin{flalign}
		&\mbox{min} \Big(c^JJ(\omega)+c^I\sum_{r\in R}I_r(\omega)+c^W\sum_{i\in S}\big(W_i(\omega)+Y_i(\omega)\big) \notag \\ 
		& + \sum_{i\in S} \mu_i(\omega) \big(a_i(\omega) - \hat{a}_i\big) + \frac{\rho }{2} \sum_{i\in S} \| a_i(\omega) - \hat{a}_i\|^2  \Big)\label{ssa}\\
		&\mbox{s.t.} \nonumber \\
		& a_j(\omega)  \geq a_i(\omega) - M (1- u_{ij}(\omega)) \qquad  \forall{i,j} \in S, j\neq{i} \label{ss1}\\
		&u_{ij}(\omega)+u_{ji}(\omega)=1\qquad   \forall{i,j} \in S, j>i  \label{ss2}\\
		&a_i(\omega)+Y_i(\omega)=B_i(\omega)\qquad  \forall i \in S \label{ss3} \\
		&a_i(\omega): integer \qquad \forall i \in S \label{ss4}\\
		&u_{ij}(\omega)\in {\{0,1\}}\qquad \forall{i,j} \in S, j\neq{i} \label{ss5}\\
		& (\ref{2d})- (\ref{2r})
\end{flalign}%
\subsection{Basic progressive hedging algorithm} 
We explain the general steps of the basic PHA, whose pseudocode is given in Algorithm \ref{alg:PH}. Letting $z$ denote the iteration counter, we initialize $z$, Lagrangian multipliers and penalty parameters in step 1. Note that $\rho^0$ represents a nonnegative constant value. In step 2, we solve SSP models for each $\omega \in \Omega$ and determine appointment time values for each patient. However, we ignore the terms with Lagrangian multipliers and penalty parameter while solving SSP models in the first iteration. We calculate consensus parameter values for each patient by taking the weighted sum of appointment time values over scenarios in step 3. Through a positive multiplier, denoted as $\alpha$, we update penalty parameter values in step 4. We update Lagrangian parameter values based on the difference between appointment time values and consensus parameter values in step 5. Finally, we check if nonanticipativity constraints are satisfied in step 6. If the constraints are satisfied, we terminate the algorithm. Otherwise, we update the iteration counter and conduct another iteration by starting from step 2.  

\begin{algorithm}
	\caption{Basic Progressive Hedging Algorithm}\label{alg:PH}
	\begin{algorithmic}[1]
		\State \textbf{Step 1: Initialize parameters:}
		\State $z=1$
		\State $\mu_i(\omega)^{(z)} = 0 \quad \forall i \in S, \forall \omega \in \Omega$
		\State ${\rho}^{(z)} = {\rho}^0$
		\State \textbf{Step 2: Solve scenario-subproblems:}
		\If{$z=1$} 
		\State Solve SSP $\forall \omega \in \Omega$ by ignoring the terms with Lagrangian multipliers and penalty parameters to determine $a_i(\omega)^{(z)} \quad\forall i \in S, \forall \omega \in \Omega $
		\Else
		\State Solve SSP $\forall \omega \in \Omega$ to determine  $a_i(\omega)^{(z)}  \quad\forall i \in S, \forall \omega \in \Omega $
		\EndIf
		\State \textbf {Step 3: Calculate consensus parameter values:} 
		\State \indent $\hat{a_{i}} = \displaystyle\sum_{\omega \in \Omega}p^\omega a_i(\omega)^{(z)} \quad \forall i \in S$
		\State \textbf {Step 4: Update penalty parameter values:}
		\If{$z>1$}  
		\State ${\rho}^{(z+1)} = \alpha {\rho}^{(z)} $
		\EndIf
		\State \textbf{Step 5: Update Lagrangian multiplier values:}
		\State\indent $\mu_i(\omega)^{(z+1)} = \mu_i(\omega)^{(z)} +{\rho}^{(z)}(a_i(\omega)^{(z)} -\hat{a_{i}}) \quad \forall i \in S, \forall \omega \in \Omega$
		\State \textbf{Step 6: Check termination criterion:}
		\If{$a_i(\omega)^{(z)} = \hat{a_{i}} \quad \forall i \in S, \forall \omega \in \Omega$}
		\State Terminate the algorithm
		\Else
		\State Set $z \leftarrow z+1$ 
		\EndIf
		\State \textbf{return} to Step 2
	\end{algorithmic}
\end{algorithm}

\subsection{Extended progressive hedging algorithm} 
Our preliminary experiments show that the basic PHA requires excessive amount of time to converge to a solution which is generally found to be of low quality. Hence, we extend the capabilities of the basic PHA by testing some enhancement ideas, and propose the algorithm that we call as \textit{extended progressive hedging algorithm (EPHA)}.     

One of the main factors that complicate the solution process is the existence of the quadratic component in the objective function of the SSP. An efficient procedure is needed to solve the SSP, as it must be solved for each scenario at each iteration of the algorithm. Therefore, we implement a linearization technique which adds cuts to each SSP independently over multiple iterations. The penalty parameter also needs a particular attention for the successful implementation of the algorithm. As it is directly represented in the objective function, and also affects Lagrangian multiplier values, the penalty parameters must be carefully updated at each iteration. Therefore, we use a dynamic penalty update method in the EPHA.  

The EPHA includes an additional step in between Step 3 and Step 4 of Algorithm \ref{alg:PH} for fixing the values of some variables. In particular, we propose a variable fixing method for the first-stage variables. Finally, we use a cycle detection algorithm to prevent cyclic behaviors. The details of the EPHA are presented in the following sections.

\subsubsection{Linearizing the objective function} 

Ignoring the penalty parameter, the quadratic component in the SSP objective function (\ref{ssa}) can be explicitly written as follows:

\begin{equation}
\label{quadexpansion}
\sum_{i\in S} \sum_{\omega\in \Omega} \| a_i (\omega) - \hat{a}_i\|^2 =  \sum_{i\in S} \sum_{\omega\in \Omega} a_i (\omega)^2 -2 \sum_{i\in S} \sum_{\omega\in \Omega} (a_i (\omega)  \hat{a}_i) + |\Omega|  \sum_{i\in S} \hat{a}_i^2
\end{equation} 
The $a_i(\omega)^2$ term in the expression in (\ref{quadexpansion}) makes the objective function quadratic. Based on our preliminary experiments, we observe that solving a quadratic mixed-integer program for each scenario at each iteration leads to a significant computational burden. Hence, we implement a linearization method which is shown to provide a good approximation of the objective function in the studies by \citet{Helseth2016} and \citet{demiretal21}.

Letting $f(x)=a_i (\omega)^2$, the method uses the first-degree of Taylor polynomial (i.e.$ \, f(x) \approx f(v) + f'(v)(x-v))$ to generate cuts as a lower bound on  $a_i (\omega)^2$ at each iteration. The quality of approximation becomes sufficient after a certain number of iterations. Since a better approximation is obtained when $v$ is closer to $x$, we set $v= a_i (\omega)^{(z-1)}$ in iteration $z$. In other words, the appointment time of patient $i$ at scenario $\omega$ found in the previous iteration is considered as an appropriate point while approximating the quadratic function on the same variable in the current iteration.

For each scenario $\omega \in \Omega $, we introduce a new variable, $h_i (\omega)$, to represent $a_i(\omega)^2$, and add the following cut to the respective SSP at iteration $z$. 

\begin{equation}
h_i(\omega)\geq (a_i(\omega)^{(z-1)})^2 + 2a_i(\omega)^{(z-1)}(a_i(\omega)-a_i(\omega)^{(z-1)}) 
\end{equation}

\subsubsection{Fixing the first-stage variables}
To reduce the computational time to solve SSPs and also accelerate the convergence process, we fix some first-stage variables before the termination criterion is satisfied. We determine the fixing strategies after conducting some preliminary experiments to observe the behaviors of the first-stage variables throughout the iterations. 

The preliminary results show that the same precedence relationships between two patients, $i, j \in S$, (i.e., $u_{ij}(\omega)$ values) may be observed in the vast majority of the SSP solutions in an early iteration. Furthermore, those $u_{ij}(\omega)$ values may not change over several iterations until the termination of the algorithm. Therefore, we fix $u_{ij}(\omega)$ values in the SSP for each $\omega \in \Omega$, when at least 80\% of the SSPs yield the same $u_{ij}(\omega)$ value for the patients $i$ and $j$ at a particular iteration. When the condition is satisfied, $u_{ij}(\omega)$ is set as the value observed in the majority of the SSP solutions. As the precedence variable fixing procedure is applied over the iterations, the range for the feasible values of the appointment time variables also gets narrower. 

Our observation on the behavior of the appointment time variables over iterations is similar to that for the precedence variables. Therefore, we fix $a_i(\omega)$ values for patient $i$ when the same value, $\tilde{a}_i$, is observed in the SSP solutions across several different scenarios. The value fixed for $a_i(\omega)$ is also selected as $\tilde{a}_i$ in that case. When the appointment time for a patient is fixed, the number of variables in the SSPs decreases, and a bound is imposed on the appointment time variables for the remaining patients. Due to these effects, the appointment time fixing procedure may significantly reduce the solution time spent for solving all SSPs in a particular iteration.        

Before implementing the appointment time fixing procedure, we set a threshold level that determines when the majority is achieved. The threshold level depends on the iteration counter. It starts from 100\% at iteration 1, and decreases by a constant amount at each iteration until iteration $limit_2$, at which point it becomes 80\%. Therefore, the amount of decrease at each iteration is set as ($20/limit_2$)\%. Next, the threshold level decreases by 10\% once at iteration $limit_3$, and again at iteration $limit_4$. Note that the high threshold levels at earlier iterations prevent the risk of the convergence to poor quality solutions. On the other hand, when the number of iterations is larger, the variability in $a_i(\omega)$ values among the SSP solutions decreases, allowing us to use lower threshold levels. The details of the variable fixing procedure are presented in Algorithm \ref{fixingalg}.

\begin{algorithm}
	\caption{Variable Fixing Procedure}\label{fixingalg}
	\begin{algorithmic}[1]
	
		\State  \textbf {if} $z=1$
		\State\indent $p^{(z)} = 100$
		\State  \textbf {end if}
		\State $\bar{u}_{ij}^{(z)} = \sum_{i\in S}\sum_{j\in S, j\neq i}u_{ij}(\omega)^{(z)}$
		\State  \textbf {for} all patient pairs  $(i,j), i\in S, j\in S, i\neq j$
		\State\indent \textbf {if} $\bar{u}_{ij}^{(z)} \geq 0.8$
		\State\indent\indent Add 	$u_{ij}(\omega) = 1  \forall \omega\in\Omega$
		\State\indent\indent Add 	$u_{ji}(\omega) = 0  \forall \omega\in\Omega$
		\State\indent  \textbf {end if}
		\State  \textbf {end for}
		\State  \textbf {for} all patient $i, i\in S$	
		\State\indent \textbf {if} $a_{ij}(\omega)^{(z)} = \tilde{a}_i$ for $p^{(z)}\%$ of the scenarios
		\State\indent\indent Add $a_i(\omega) = \tilde{a}_i  \forall \omega\in\Omega$
		\State\indent  \textbf {end if}
		\State  \textbf {end for}
		\State  \textbf {if} $z\leq limit_2$
		\State\indent $p^{(z+1)} = p^{(z)}- 20/limit_2$
		\State  \textbf {else if} $z =limit_3$
		\State\indent $p^{(z+1)} = p^{(z)}-10$
		\State  \textbf {else if} $z =limit_4$
		\State\indent $p^{(z+1)} = p^{(z)}-10$
		\State  \textbf {end if}
	\end{algorithmic}
\end{algorithm}

During the implementation of the basic PHA, the appointment time values of all patients in the SSPs for all scenarios may stay the same over consecutive iterations. To prevent such a cyclic behavior and the risk of non-convergence, we implement the cycle detection approach proposed by \citet{Watson2011}. The approach detects cycles by tracking the Lagrangian multipliers for each patient, whose values change even if a cycle occurs. When a cycle is detected, we fix the appointment time value to the $\hat{a}_i$ value found at that iteration. Note that the cycle detection is implemented at each iteration after iteration $limit_2$.

The cycle detection and variable fixing procedures we implement ensure the convergence of the EPHA. However, even though it occurs rarely, the EPHA may run through an extensive number of iterations in some instances, because the appointment times of a number of patients can not be fixed. We observe in our preliminary experiments that after a certain iteration (we call it as $limit_5$), the unfixed appointment times change slightly. Therefore, to prevent running the algorithm for several minutes only to gain a minor improvement in the objective value, we review the unfixed patient appointment times every $controliter$ iterations after iteration $limit_5$. If the number of unfixed appointment times does not change in the last $controliter$ iterations, we fix the appointment time of each patient $i$ whose appointment time was not fixed. Note that the fixed value is chosen as the $\hat{a}_i$ value in that iteration.    

\subsubsection{Penalty Parameter Update Method}
The selection of penalty parameter values is critical to ensure that the EPHA converges to a good quality solution in a reasonable amount of time. Setting very high values for $\rho$ enforces the convergence of appointment time values in the SSPs for different scenarios to the related consensus parameter values in the early iterations. However, the consensus parameter values are estimated progressively, and their estimations can be of low quality in earlier iterations. Note that the initial estimation of consensus parameters are obtained by solving SSPs without any non-anticipativity constraints. Therefore, very high values of $\rho$ may result in low quality model solutions. On the other hand, very low $\rho$ values would delay the convergence of the appointment time values, and hence increase the number of iterations. As the number of iterations increases, the appointment time values can be better estimated at the expense of increased computational time. Hence, the EPHA can converge to a good quality solution after extensive number of iterations. Therefore, the value of $\rho$ must be set carefully at each iteration by examining its relationship with the solution quality and run time. 

The previous research shows that using a dynamic update method to change parameter values provide significant advantages in terms of solution quality and run time \citep{Hvattum2009, Gul2015, demiretal21}. We enhance the penalty update method suggested by \citet{demiretal21}, which is an extended version of \citet{Hvattum2009}. 

The method uses two parameters that help examine the convergence pattern of both primal and dual solutions. The first parameter, denoted by $\Delta_d$, represents the squared sum of the differences between appointment time values in the SSP solution for each scenario and the relevant consensus parameter. This parameter provides information on the convergence level at a particular iteration. If $\Delta_d$ increases from one iteration to another, this means the SSP solutions are moving away from the consensus condition. The second parameter, denoted by $\Delta_p$, considers the squared sum of the difference between the consensus parameter value in the current iteration and that in the previous iteration, and sum the differences over all patients. An increase in $\Delta_p$ implies that the EPHA is moving from the previous target to another target of convergence. 

The method first checks the change in $\Delta_d$ value with respect to that of the previous iteration. If the change is positive, $\rho$ is multiplied by $\alpha>1$ to ease convergence. Otherwise, it implies that the EPHA solution is approaching a consensus point. Next, the method checks the change in $\Delta_p$ value from the previous iteration to the current one. If the change is positive, $\rho$ is multiplied by $\frac{1}{\alpha}$ to prevent the convergence to a poor solution. Otherwise, the method keeps $\rho$ stable, since the decrease in  $\Delta_p$ value means that the oscillation from one consensus point to another one is fading. 

We impose a three-phase limit on $\rho$ while implementing the penalty update method to reduce the risk of immature termination of the EPHA, which may occur due to the following reasons: (i) the number of cuts added to approximate the quadratic term becomes sufficient only after a certain number of iterations, or (ii) a rapid increase in $\rho$ may activate variable fixing procedures unnecessarily. Until iteration $limit_1$, $\rho$ can not exceed $\rho^{u1}$. After that point, the upper limit is updated to $\rho^{u2}$, which is greater than $\rho^{u1}$. The second upper limit is valid until iteration $limit_5$. Afterwards, the upper limit is updated to $\rho^{u3}$, a value greater than $\rho^{u2}$. Note that $limit_5$ is the threshold level also used within the variable fixing procedure. The details of the penalty parameter update method is shown in Algorithm \ref{penupdate}.

\begin{algorithm}
	\caption{Penalty Update Method with Three Upper Limits}\label{penupdate}
	\begin{algorithmic}[1]
		
		\State  $\Delta_d^{(z)} =\displaystyle\sum_{i \in S} \displaystyle\sum_{\omega \in \Omega} (a_i(\omega)^{(z)} - \hat {a}_i^{(z)})^2$
		\State  $\Delta_p^{(z)} =\displaystyle\sum_{i \in S} (\hat {a}_i^{(z)} - \hat {a}_i^{(z-1)})^2$
		\State \textbf{if} $z \leq limit_1$
		\State\indent $\rho^u =\rho^{u1}$
		\State \textbf{else if} $z \leq limit_5$
		\State\indent $\rho^u =\rho^{u2}$
		\State \textbf{else} 
		\State\indent $\rho^u =\rho^{u3}$
		\State \textbf{if} $\Delta_d^{(z)} - \Delta_d^{(z-1)} > 0 \quad \&  \quad \rho^{(z)} < \rho^u $
		\State\indent $\rho^{(z+1)} = \alpha \rho^{(z)}$
		\State \textbf{else if} $\Delta_d^{(z)} - \Delta_d^{(z-1)} > 0$
		\State\indent $\rho^{(z+1)} = \rho^u$
		\State \textbf{else if} $\Delta_p^{(z)} - \Delta_p^{(z-1)} > 0$ 
		\State\indent $\rho^{(z+1)} = \frac{1}{\alpha \rho^{(z)}}$
		\State \textbf{else if} $ \rho^{(z)} > \rho^u $
		\State\indent $\rho^{(z+1)} = \rho^u$
		\State \textbf{else}
		\State\indent $\rho^{(z+1)} = \rho^{(z)}$
		
	\end{algorithmic}
\end{algorithm}

\section{Experimental Study}
We test the EPHA using the data set of an outpatient procedure center in Mayo Clinic \citep{Gul2011}. We create problem instances based on the surgeries performed at the urology department. We sample induction and surgery durations from the data set which includes the records of 1963 patients whose surgeries were performed over the data collection period. Urology surgeries are grouped into five acuity levels in the data set. The number of surgeries of each acuity level and descriptive statistics related to induction and surgery durations are given in Table \ref{urodata}. We sample induction and surgery durations independently for each acuity level. We generate turnover durations by assuming a uniform distribution with lower and upper limits as 15 and 30 minutes, respectively. We implement the algorithms in Microsoft Visual C++ 2019 using CPLEX 12.8 Concert Technology. We conduct the experiments on an Intel(R) Xeon(R) E-2246G computer with six-core processor running at 3.60 GHz and 16GB RAM. 

Each experiment is conducted on an instance set that consists of 10 instances. An instance set represents the case where 7 patients are scheduled for 2 IRs and 3 ORs. We ensure that occurrence frequencies of acuity levels in a daily patient list are proportional to the patient count values for each acuity level in the data set given in Table \ref{urodata}. The instances in a given instance set differ from each other based on the induction and surgery durations. In other words, different set of induction and surgery duration values are sampled for a patient for each instance. The number of duration values sampled for each patient in an instance is set to 50 (i.e. $|\Omega|=50$) in all experiments, except where we assess the optimality gap.

We consider the trade-off between idle time and waiting time by setting different values for the cost coefficients, $c^I, c^J, c^W$. Instead of setting actual cost values for the coefficients, we use them as trade-off parameters whose values add up to 1 in each experiment.     

We set the values of the parameters used within the variable fixing and penalty update methods based on our preliminary experiments as $limit_1=25$, $limit_2=50$, $limit_3=60$, $limit_4=70$ and $limit_5=90$, $controliter=100$.
\begin{table}
	\centering
	\caption{Number of patients, mean and standard deviations of induction and surgery durations (in minutes) for each acuity level}
	\label{urodata}
	\small{
		\begin{tabular}{|c|c|c|c|c|c|}
		
			\cline{3-6}    \multicolumn{2}{r|}{} & \multicolumn{2}{c|}{\textbf{Induction}} & \multicolumn{2}{c|}{\textbf{Surgery}} \\
			\hline
			\textbf{Acuity level} & \textbf{Count}  & \textbf{Mean} & \textbf{Std. deviation}  & \textbf{Mean} & \textbf{Std. deviation}  \\
			\hline
			1 & 329  & 23.65  & 5.66 & 29.65 & 20.52 \\
		
			2 & 640 & 13.82 & 6.07 & 17.48 & 8.58 \\
			
			3 & 153 & 29.03 & 6.78 & 109.12 & 42.96 \\
		
			4 & 345 & 24.97 & 5.67 & 30.88 & 14.12 \\
		
			5 & 496 & 28.20 & 7.42 & 52.12 & 32.56 \\
			\hline
	\end{tabular}}
\end{table}%

\subsection{Benefit of variable fixing and penalty update methods}

We assess the value of implementing variable fixing procedure and penalty update method. We set $c^I, c^J, c^W$ as 0.5, 0.25, 0.25 in these experiments, respectively. Table \ref{varfix} compares the performance of the algorithm with and without variable fixing procedure over 10 instances. The results clearly illustrate the significant benefit of the variable fixing procedure. The table shows that when the variable fixing mechanism is implemented, the objective value improves by 6.1\% on average. Furthermore, the algorithm terminates due to time limit of 3 hours if the variable fixing procedure is not incorporated into the algorithm. On the other hand, the EPHA finds the solutions in around 14 minutes when the variable fixing procedure is implemented.   

\begin{table}
	\centering
	\caption{Benefit of variable fixing procedure in terms of objective value and run time (in seconds)}
	\label{varfix}
	\small{
		\begin{tabular}{|c|c|c|c|c|}
			
			\cline{2-5}    \multicolumn{1}{r|}{} & \multicolumn{2}{c|}{\textbf{Objective Value}} & \multicolumn{2}{c|}{\textbf{Run Time}} \\ \hline
			\textbf{Instance \#}  & \textbf{Fixing} & \textbf{No Fixing} & \textbf{Fixing} & \textbf{No Fixing} \\ \hline
			1  &    73.27          & 79.58    & 796.20   &  10800.00           \\ 
			2  &    71.76          & 75.40   & 1153.85   &  10800.00           \\ 
			3  &    70.77          & 78.37    & 715.11   &  10800.00           \\ 
			4  &    67.36          & 67.44    & 645.28   &  10800.00           \\ 
			5  &    75.91          & 75.87   & 678.02   &  10800.00           \\ 
			6  &    77.40          & 78.59    & 852.65   &  10800.00           \\ 
			7  &    70.41          & 71.28   & 891.20   &  10800.00           \\ 
		    8  &    76.27          & 93.76  & 715.64   &  10800.00           \\ 
			9  &    72.25          & 77.82    & 923.53   &  10800.00           \\ 
			10 &    75.03          & 79.86    & 893.79   &  10800.00           \\ \hline
			\textbf{Average} & 73.04 & 77.79  & 826.53   &  10800.00           \\ \hline
	\end{tabular}}
\end{table}

We next show the benefit of implementing the penalty update method within the EPHA. As mentioned in Section 4, our penalty update method considers 3 upper limits ($\rho^{u1},\rho^{u2}, \rho^{u3}$). To illustrate the benefit of using multiple upper limits, we first conduct experiments using the method with 1 ($\rho^{u1}$) and 2 upper limits ($\rho^{u1},\rho^{u2}$) based on the baseline instances. Table \ref{penup1} shows that when the number of imposed upper limits increases from 1 to 2, the run time improves significantly (by 27\%), while the objective value improves by 0.08\%. We then conduct tests using 2  ($\rho^{u1},\rho^{u2}$) and 3 upper limits ($\rho^{u1},\rho^{u2}, \rho^{u3}$) on another instance set selected to ensure that the EPHA iteration counter exceeds $limit_5$ such that the limit $\rho^{u3}$ can be imposed in all instances. Table \ref{penup2} illustrates that increasing the number of upper limits from 2 to 3 provides a greater improvement in terms of the objective value (by 1.76\%). Since the run time also improves by 26\%, we continue with using 3 upper limits while implementing the penalty update method in the remainder of the experiments.           

\begin{table}
	\centering
	\caption{Benefit of the penalty update method (using 1 or 2 upper limits) in terms of objective value and run time (in seconds)}
	\label{penup1}
	\small{
		\begin{tabular}{|c|c|c|c|c|}
			
			\cline{2-5}    \multicolumn{1}{r|}{} & \multicolumn{2}{c|}{\textbf{Objective Value}} & \multicolumn{2}{c|}{\textbf{Run Time}} \\ \hline
			\textbf{Instance \#}  & \textbf{2} & \textbf{1} & \textbf{2} & \textbf{1} \\ \hline
			1  &    73.27          & 73.57    & 796.20   & 2894.29           \\ 
			2  &    71.76          & 71.73   & 1153.85   &  1443.02           \\ 
			3  &    71.77          & 71.44    & 715.11   &  1024.66           \\ 
			4  &    67.36          & 67.46    & 645.28   &  955.09           \\ 
			5  &    75.91          & 75.85   & 678.02   &  931.84           \\ 
			6  &    77.40          & 77.45    & 852.6   & 929.56           \\ 
			7  &    70.41          & 70.33   & 891.20   &  1020.38           \\ 
		    8  &    76.27          & 76.53  & 715.64   &  1811.83           \\ 
			9  &    72.25          & 72.21    & 923.53   &  885.00           \\ 
			10 &    75.03          & 74.45    & 893.79   &  1169.84           \\ \hline
			\textbf{Average} & 73.04 & 73.10  & 826.53   &  1306.55           \\ \hline
	\end{tabular}}
\end{table}

\begin{table}
	\centering
	\caption{Benefit of the penalty update method (using 2 or 3 upper limits) in terms of objective value and run time (in seconds)}
	\label{penup2}
	\small{
		\begin{tabular}{|c|c|c|c|c|}
			
			\cline{2-5}    \multicolumn{1}{r|}{} & \multicolumn{2}{c|}{\textbf{Objective Value}} & \multicolumn{2}{c|}{\textbf{Run Time}} \\ \hline
			\textbf{Instance \#}  & \textbf{3} & \textbf{2}  &\textbf{3} & \textbf{2}  \\ \hline
			1  &   78.05           &  78.11   & 1380.02   &  1488.61          \\ 
			2  &   80.97           & 80.99  & 1390.41   &  1677.11         \\ 
			3  &   79.07           &  79.19   & 1266.22   &  1433.48          \\ 
			4  &   76.57           &  76.34   & 1556.31   &  2367.65          \\ 
			5  &   81.98           &  81.98   & 1674.81   &  1531.63        \\ 
			6  &   78.16           &  78.53   & 1599.59   &  2327.15          \\ 
			7  &   79.94           &  79.94   & 1614.65   &  1619.08         \\ 
			8  &   84.61           &  98.52   & 3324.43   &  6171.85         \\ 
			9  &   80.99           &  81.23   & 1674.25   &  2072.15        \\ 
			10 &   72.86           &  72.56   & 1678.84   &  2630.24         \\ \hline
			\textbf{Average} & 79.32 & 80.74  &  1715.95  &  2331.88         \\ \hline
	\end{tabular}}
\end{table}

\subsection{Comparison of the EPHA results with optimal solutions}
We assess the performance of the EPHA by comparing the algorithm with CPLEX on three different instance sets. The sets are differentiated from each other by varying cost values (i.e., the number of surgeries from each acuity level does not change). Since resource utilization is given higher priority over patient satisfaction by the decision makers in hospitals, we ensure that the unit cost of waiting time does not exceed the total unit costs of OR and IR idle times in these experiments. Furthermore, we consider only the realistic cases where the unit cost of IR idle time does not exceed that of OR idle time. Note that all other cases are considered in Section 5.4. In particular, we consider the following values for the triplet, $c^I, c^J, c^W$: (0.33, 0.33, 0.34),(0.5, 0.25, 0.25), (0.65, 0.10, 0.25). We set $|\Omega|=10$ in these experiments, as CPLEX cannot find optimal solutions for larger values of scenario sizes.

Table \ref{ephagap} shows that the average gap between the EPHA and CPLEX over all instances is 2.39\%, while it can be as low as 1.05\% for an instance set. The EPHA finds the solution in about 1.5 minutes on average, while the CPLEX requires more than 41 minutes to obtain the optimal solution.  

\begin{table}
	\centering
	\caption{Comparison of the EPHA with CPLEX with respect to the objective values and run times (in seconds) for each instance set}
	\label{ephagap}
	\small{
	\begin{tabular}{|c|c|c|c|c|c|c|}
	
		\cline{3-7}    \multicolumn{2}{r|}{} & \multicolumn{3}{c|}{\textbf{Objective Value}} & \multicolumn{2}{c|}{\textbf{Run Time}} \\ \hline
		\textbf{Set \#} & \textbf{$c^I, c^J, c^w$} & \textbf{EPHA} & \textbf{CPLEX} & \textbf{Difference (\%)}& \textbf{EPHA} & \textbf{CPLEX} \\ \hline
		1 & 0.33, 0.33, 0.34 &          59.05      & 56.2   & 4.73 &  123.86  &   2317.97          \\ 
		2 & 0.5, 0.25, 0.25  &      67.27         &      66.55         & 1.05& 71.12  & 3236.03       \\ 
		3 & 0.65, 0.10, 0.25  &          73.92      &     72.91       & 1.39 &  88.47  & 1844.64        \\ 
	    Average &   &      66.74          &    65.22        & 2.39 &94.48&   2466.19          \\ \hline
	\end{tabular}}
\end{table}

\subsection{Comparison of the EPHA with scheduling heuristics}
We compare the EPHA with the commonly used scheduling heuristics in the earlier literature on surgery or chemotherapy appointment scheduling \citep{Gul2011, Gul2018,demiretal21}. In particular, we make comparisons with combinations of sequencing and job hedging heuristics.

 Using sequencing heuristics, we sort the surgeries based on one of the three rules: (1) smallest  expected induction duration  first (SPT), (2) largest expected induction duration first (LPT), and (3) smallest variance of induction first (VAR). Note that we choose induction duration as the criterion, since the sequencing variables in the two-stage SMIP model exist only for the induction activity. Furthermore, as can be verified by Table \ref{urodata}, SPT and LPT sequences do not change even if the patients are sequenced based on surgery duration or total duration. The sequence for VAR rule may change in such a case, but our preliminary experiments show that the induction duration is a better criterion.  
 
 After sequencing the surgeries, we calculate the appointment times of patients. Initially, the appointment times of the first patients of each IR are set as 0 (i.e., the beginning of the shift). We then identify the patient whose induction is expected to finish the earliest among all patients in the IR based on their induction durations, which are estimated based on the job hedging procedure. Let that patient be indexed by $i$. Next, the first patient among the non-treated patients is identified (let index $j$ denote that patient). The appointment time of patient $j$ is then set by adding the estimated induction duration of patient $i$ to the appointment time of patient $i$. In other words, the appointment time of patient $j$ is set equal to the planned induction finish time of patient $i$. Due to the trade-off between waiting time and idle time, the estimation of induction duration is critical. We estimate the induction duration for each acuity level separately, while applying job hedging. Job hedging procedure allows to use different percentiles of the relevant distribution while estimating the induction duration. In our experiments, we test $50^{th}$ (median), $60^{th}$, $70^{th}$, $80^{th}$, and $90^{th}$ percentiles. Instead of fitting general distributions, we sort the induction durations for each acuity level in our data set, and identify the values corresponding to the tested percentiles based on the sample. When, the appointment time of patient $j$ is set, the next patient among the non-treated patients is identified similarly. We set the appointment time of this patient, and then follow the same procedure, iteratively, to set the appointment times of all patients in the daily list.    
 
 \begin{table}
 	\centering
 	\caption{Comparison of EPHA with scheduling heuristics}
 	\label{jobhedging}
 	\small{
 		\begin{tabular}{|c|c|c|c|}
 			\cline{2-4}    \multicolumn{1}{r|}{}  &\multicolumn{3}{c|}{\textbf{Difference (\%)}} \\ \hline
 			\textbf{Percentile}  &\textbf{SPT} & \textbf{LPT} & \textbf{VAR}  \\ \hline
 			50\%   &42.33 &   71.93 &    54.38               \\ 
 			 60\%  &40.29  &    70.66         &    54.25                 \\ 
 			 70\%     &39.20   &  70.22           &    53.76                 \\ 
 			  80\%     & 37.89 &   70.73          &      52.04              \\ 
 			  90\%     &  37.52 &   66.72         &    54.76               \\ 
 			  \hline
 			 \textbf{ Average}    & 39.45  &  70.05          &     53.84                \\
 			 \hline
 	\end{tabular}}
 \end{table}
 
We set $c^I, c^J, c^W$ as 0.5, 0.25, 0.25 in these experiments, respectively. The average EPHA objective value over 10 instances in the set is found as 73.04. Table \ref{jobhedging} provides the percentage improvement that the EPHA solution provides over the heuristic solutions with respect to the objective value. Among the sequencing heuristics, SPT performs the best on average. However, the EPHA improves SPT solutions by almost 40\%. The improvement over VAR and LPT solutions are around 54\% and 70\%, respectively. Note that increasing job hedging level results in better results when SPT is used to sequence patients, but we cannot generalize this observation to other sequencing rules. The overall results indicate that it is highly essential to use a sophisticated algorithm such as the EPHA to schedule surgeries for surgical suites that follow parallel processing principle.
 
\subsection{Sensitivity analysis on the model and problem parameters}

In this section, we present our experiments on analyzing the sensitivity of the results of EPHA with respect to changing unit costs of the three objectives and the number of induction rooms.

\subsubsection{Impact of the unit costs of the three objectives}

Our first sensitivity analysis covers the case where we vary the relative unit costs of OR idle time, IR idle time and patient waiting time to observe the trade-off between these three objectives. For this end, we use 12 different combinations of $c^I$, $c^J$ and $c^W$ by setting two of the costs equal to each other and the remaining third to be 10\%, 50\%, 200\%, and 1000\% of the former two. We then normalize these unit costs to ensure that they add up to 1.

\begin{table}
    \centering
    \caption{Average objective values, run times (in seconds), idle times and patient waiting times (in minutes) of EPHA for different unit costs of the three objectives}
    \label{tab:lambda}\small
    \begin{tabular}{|c|c|c|c|c|c|c|c|}
        \hline
        $c^I$ & $c^J$ & $c^W$ & \textbf{Objective} & \textbf{Run time} & \textbf{OR idle time} & \textbf{IR idle time} & \textbf{Waiting time} \\
        \hline
        0.833	&	0.083	&	0.083	&	78.14	&	626.46	&	79.17	&	63.41	&	82.60	\\
        0.500	&	0.250	&	0.250	&	72.51	&	907.00	&	92.17	&	49.25	&	56.47	\\
        0.476	&	0.476	&	0.048	&	62.90	&	699.92	&	86.37	&	35.65	&	100.68	\\
        0.476	&	0.048	&	0.476	&	69.27	&	807.75	&	100.47	&	61.52	&	38.83	\\
        0.400	&	0.400	&	0.200	&	64.39	&	1003.45	&	98.61	&	37.93	&	48.87	\\
        0.400	&	0.200	&	0.400	&	65.58	&	1131.80	&	102.45	&	44.33	&	39.33	\\
        0.250	&	0.500	&	0.250	&	56.96	&	1023.33	&	128.20	&	25.41	&	48.82	\\
        0.250	&	0.250	&	0.500	&	58.08	&	1240.47	&	130.73	&	44.05	&	28.78	\\
        0.200	&	0.400	&	0.400	&	53.04	&	1189.69	&	134.61	&	32.07	&	33.21	\\
        0.083	&	0.833	&	0.083	&	35.84	&	876.44	&	228.88	&	13.15	&	69.76	\\
        0.083	&	0.083	&	0.833	&	28.92	&	1334.53	&	178.77	&	65.73	&	10.25	\\
        0.048	&	0.476	&	0.476	&	31.31	&	1134.13	&	162.15	&	24.51	&	25.03	\\
        \hline															
    \end{tabular}
\end{table}

For each of the 12 unit cost combinations, Table \ref{tab:lambda} provides the average objective function value, computational time, OR idle time, IR idle time, and patient waiting time for each instance set. Focusing on each of the last three columns yields insights into the sensitivity of each of these objectives to the change in unit costs. Among these, OR idle time ranges from 79.17 in the case where its unit cost is 10 times that of the other two objectives to 228.88 minutes when the cost per minute of IR idle time is 10 times that of OR idle time and patient waiting time. Similarly, IR idle time ranges from 13.15 to 63.41 minutes, and patient waiting time varies from 10.25 to 100.68 minutes. Whereas the range is largest for OR idle time, patient waiting time sees the highest relative increase from average, with a value twice as long in the worst case. The variations of all three objectives with regard to unit cost changes underlines the significance of calculating these costs accurately when determining the surgery schedules.

When a pairwise analysis of changes is made between the objectives, one observes a negative correlation between each pair, which shows that these three objectives are in total conflict with one another, pointing to the importance of incorporating all three objectives into the model. The most significant trade-off occurs between OR and IR idle time, with a correlation coefficient of -0.44. The corresponding coefficients for (i) OR idle time and waiting time, and (ii) IR idle time and waiting time are -0.30 and -0.28, respectively.

\subsubsection{Impact of the number of induction rooms}

An important decision for surgery scheduling in the existence of parallel processing is the number of IRs in the system. While a large number of IRs may improve the performance measures, it has been observed in the literature that using more IRs leads to the need to hire more anesthesia and nursing personnel \citep{Harders2006}. If the management decides to staff additional IRs with the existing personnel, this would induce opportunity cost, as the personnel could be utilized elsewhere in the hospital. Hence, care should be taken into consideration in determining the number of IRs. To analyze the effect of the number of IRs on system performance, we extend our experiments to involve (i) a single IR and (ii) three IRs, in addition to the two-IR case in our baseline setting.

\begin{table}
    \caption{Average objective values, run times (in seconds), idle times and patient waiting times (in minutes) of EPHA for varying number of induction rooms}
    \label{tab:irs}
    \centering \small
    \begin{tabular}{|c|c|c|c|c|c|} 
        \hline
        \textbf{Number of IRs}	&	\textbf{Objective}	&	\textbf{Run Time}	&	\textbf{OR Idle Time}	&	\textbf{IR Idle Time}	&	\textbf{Waiting Time}	\\	\hline
        1	&	129.14	&	802.21	&	187.42	&	53.49	&	88.23	\\	
        2	&	72.51	&	907.00	&	92.17	&	49.25	&	56.47	\\	
        3	&	67.68	&	1707.79	&	83.42	&	55.19	&	48.70	\\	\hline
    \end{tabular}
\end{table}

Table \ref{tab:irs} summarizes the results with varying numbers of IRs. The results show that adding a second IR decreases the overall objective by around 45\%, whereas the benefit of adding the third IR is less than 7\%. As expected, the law of diminishing returns applies in this case. We observe similar outcomes when comparing the OR idle times and patient waiting times for the three settings, where the reduction is more substantial with the addition of the second IR, but decreases significantly with the third one. IR idle time does not change substantially with increasing number of IRs. While this may not be an intuitive result, it is due to the way idle time is calculated based on the IR closure times. As more IRs are added, they become less utilised. However, at the same time, closure times are shortened, offsetting the increase in the idle time and resulting in similar numbers for all settings.

\subsection{Comparison of parallel processing and serial processing}

As opposed to a serial system where both induction and turnover are performed within the OR, the parallel processing approach allows performing these activities simultaneously in different rooms, increasing the utilization of ORs. The main advantage with a serial system would be a decrease in patient waiting times, as a parallel system involves both waiting for the IR and the OR, whereas a serial system only includes waiting for the OR. On the other hand, in a serial system, the need for carrying out induction and turnovers in the same OR implies that the room will be closed at a later time compared to a parallel processing system.

\begin{table}
    \caption{Patient waiting times and OR closure times (in minutes) of EPHA under different parallel and serial processing settings}
    \label{tab:parvsseq}
    \centering \footnotesize
    \begin{tabular}{|c|c|c|c|c|c|c|c|c|c|}
    \hline
        & \multicolumn{3}{c|}{\textbf{2IR-parallel}} & \multicolumn{3}{c|}{\textbf{3IR-parallel}} & \multicolumn{3}{c|}{\textbf{Serial}} \\
        \hline
        $q_i(\omega)/e_i(\omega)$ & 0.5 & 1 & 2 & 0.5 & 1 & 2 & 0.5 & 1 & 2 \\
        \hline
        \textbf{Waiting time}	&	51.85	&	58.05	&	67.34	&	47.74	&	53.33	&	67.49	&	36.12	&	39.65	&	52.78	\\
        \textbf{OR1 closure}	&	92.31	&	107.36	&	144.61	&	85.68	&	103.33	&	141.48	&	109.07	&	127.49	&	167.30	\\
        \textbf{OR2 closure}	&	153.54	&	178.87	&	237.19	&	143.46	&	168.66	&	227.59	&	166.85	&	193.78	&	258.12	\\
        \textbf{OR3 closure}	&	252.02	&	285.41	&	361.02	&	246.65	&	286.57	&	349.74	&	301.78	&	340.29	&	413.22	\\
        \hline																			
    \end{tabular}
\end{table}

To compare the parallel processing approach against an equivalent serial system, we consider two parallel systems, one with two and one with three IRs, as well as a serial system without an IR. In all three cases, there are three ORs. We also vary the ratio between the turnover and induction time as 0.5, 1, and 2 to observe whether the conclusions change with regard to the durations of these activities. As the serial system does not involve any IRs, IR idle time cannot be used as a performance measure. Due to the same reason, OR idle time ceases to be a fair comparison measure. Consequently, we use the closure times of the three ORs, in addition to the total patient waiting time, to compare these two systems. The results of our experiments are summarized in Table \ref{tab:parvsseq}.

Table \ref{tab:parvsseq} shows the trade-off between patient waiting time and OR closure times between the two systems quite clearly. The serial system yields up to 20 and 15 minutes of total patient waiting time savings over a 2IR and a 3IR system, respectively. On the other hand, these savings are more than offset in the parallel processing system by the savings in OR closure times. With two IRs, the total savings varies between 80 and 95 minutes. Setting up three IRs increases the savings to between 100 and 120 minutes. The findings are robust over varying ratios of the ratio of turnover to induction time. These findings underline the increase in overall efficiency when using a parallel processing system.

\subsection{Estimation of the Value of Stochastic Solution}

Our last set of experiments focus on estimating the value of stochastic solution (VSS). VSS aims to assess the additional benefits obtained by considering the uncertainty of the problem environment when determining the decisions in the first stage of the problem. To estimate the VSS, we first solve the mean value problem. This is obtained by setting all the uncertain parameters at their expected values and solving a deterministic model. The objective values arising from the first-stage decisions of the mean value problem are compared against those of the EPHA to calculate the relative VSS values, given in Table \ref{tab:VSS}.

\begin{table}
    \centering
    \small
    \caption{Average objective values for the mean value solution (MV), EPHA solution, and relative VSS for the instances}
    \label{tab:VSS}
    \begin{tabular}{|c|c|c|c|}
        \hline
        \textbf{Instance}	&	\textbf{MV}	&	\textbf{EPHA}	&	\textbf{\%VSS}	\\	\hline
        1	&	78.40	&	74.20	&	5.67\%	\\	
        2	&	77.74	&	73.71	&	5.46\%	\\	
        3	&	82.92	&	76.06	&	9.03\%	\\	
        4	&	75.76	&	70.17	&	7.97\%	\\	
        5	&	75.64	&	70.46	&	7.36\%	\\	
        6	&	77.83	&	74.31	&	4.74\%	\\	
        7	&	77.38	&	73.77	&	4.89\%	\\	
        8	&	77.19	&	70.83	&	8.97\%	\\	
        9	&	84.97	&	80.21	&	5.93\%	\\	
        10	&	73.24	&	69.11	&	5.98\%	\\	
        \hline
    \end{tabular}
\end{table}

For all of our 10 instances generated based on the baseline parameter settings, the relative VSS ranges between 5\% and 10\%. While this may not appear to be a considerable change, it points to a substantial savings in the OR costs. The findings from Table \ref{tab:VSS} underline the importance of considering the uncertainty in the surgery, induction and turnover times when deciding on the appointment times of the patients for a surgical suite following parallel processing principles.

\section{Conclusions}
In this study, we propose a two-stage SMIP model for the problem of scheduling surgeries for multiple ORs and IRs that function based on parallel processing principle. The parallel processes in our problem refer to the concurrent implementation of induction in an IR and turnover in an OR. We consider the uncertainty in induction, surgery and turnover durations while making scheduling decisions. Using the two-stage SMIP model, we first sequence patients and set appointment times for surgeries. After the realization of uncertain durations, we assign patients to IRs independently for each scenario. We show that the special structure of the model allows the use of an optimal myopic policy for IR assignment decisions. We minimize the expected total cost of patient waiting time, OR idle time and IR idle time, whose values are calculated at the second stage of the two-stage SMIP model. We add valid inequalities by imposing bounds on variables and formulating symmetry-breaking constraints to enhance the model formulation.

We implement an extended version of the PHA, which we call as EPHA, to solve the two-stage SMIP model. We propose a penalty update method and a variable fixing mechanism within the EPHA. Besides, we add cuts to linearize the quadratic objective function of the scenario subproblems and implement a cycle detection mechanism using the procedures proposed in the earlier literature. 

We illustrate the advantages provided by the EPHA over the basic PHA by conducting a comprehensive set of experiments using the data set of a large academic hospital in the US. We compare the EPHA with CPLEX and show that the algorithm provides near-optimal solutions within a reasonable amount of time. We also compare the EPHA with several combinations of simple sequencing and job hedging heuristics and indicate that a sophisticated solution approach is necessary to solve the two-stage SMIP model. We assess the VSS to emphasize the importance of considering uncertainty in induction, surgery and turnover durations. 

Based on the results obtained through near-optimal schedules, we show that a parallel processing system significantly outperforms the serial system in terms of OR closure times. The finding is consistent with the results of the actual trials conducted in the hospitals without optimizing surgery schedules. The lower OR closure times in parallel processing systems allows a surgical suite manager to increase daily throughput. On the other hand, the average total patient waiting time for a parallel processing system is found to be worse than that of the serial system. Therefore, the manager must carefully consider the trade-off between patient satisfaction and daily throughput while choosing a particular setting for the flow of patients. As the number of IRs increases, the manager can better benefit from the parallel processing principle. However, the additional benefits must be weighed against the additional staffing costs while determining the number of IRs.

In this study, we make the assumption that the patient-to-OR assignment decisions are made in advance. We plan to relax this assumption and create a model for a more flexible case in a future study.       

\bibliographystyle{informs2014}
\bibliography{references}

\begin{thebibliography}{38}
\providecommand{\natexlab}[1]{#1}
\providecommand{\url}[1]{\texttt{#1}}
\providecommand{\urlprefix}{URL }

\bibitem[{Ahmadi-Javid et~al.(2017)Ahmadi-Javid, Jalali, \protect\BIBand{}
  Klassen}]{AhmadiJavid2017}
Ahmadi-Javid A, Jalali Z, Klassen KJ (2017) Outpatient appointment systems in
  healthcare: A review of optimization studies. \emph{European Journal of
  Operational Research} 258(1):3--34.

\bibitem[{Atighehchian et~al.(2020)Atighehchian, Sepehri, Shadpour,
  \protect\BIBand{} Kianfar}]{Atighehchian2020}
Atighehchian A, Sepehri MM, Shadpour P, Kianfar K (2020) A two-step stochastic
  approach for operating rooms scheduling in multi-resource environment.
  \emph{Annals of Operations Research} 292:191--214.

\bibitem[{Bai et~al.(2017)Bai, Storer, \protect\BIBand{} Tonkay}]{Bai2017}
Bai M, Storer R, Tonkay G (2017) A sample gradient-based algorithm for a
  multiple-or and pacu surgery scheduling problem. \emph{IISE Transactions}
  49(4):367--380.

\bibitem[{Batun et~al.(2011)Batun, Denton, Huschka, \protect\BIBand{}
  Schaefer}]{Batun2011}
Batun S, Denton BT, Huschka TR, Schaefer AJ (2011) Operating room pooling and
  parallel surgery processing under uncertainty. \emph{INFORMS journal on
  Computing} 23(2):220--237.

\bibitem[{Birge \protect\BIBand{} Louveaux(2011)}]{birge2011introduction}
Birge JR, Louveaux F (2011) \emph{Introduction to stochastic programming}
  (Springer Science \& Business Media).

\bibitem[{Cayirli \protect\BIBand{} Veral(2003)}]{Cayirli2003}
Cayirli T, Veral E (2003) Outpatient scheduling in health care: a review of
  literature. \emph{Production and Operations Management} 12(4):519--549.

\bibitem[{Crainic et~al.(2011)Crainic, Fu, Gendreau, Rei, \protect\BIBand{}
  Wallace}]{Crainic2011}
Crainic TG, Fu X, Gendreau M, Rei W, Wallace SW (2011) Progressive
  hedging-based metaheuristics for stochastic network design. \emph{Networks}
  58(2):114--124.

\bibitem[{Demir et~al.(2021)Demir, Gul, \protect\BIBand{}
  {\c{C}}elik}]{demiretal21}
Demir NB, Gul S, {\c{C}}elik M (2021) A stochastic programming approach for
  chemotherapy appointment scheduling. \emph{Naval Research Logistics (NRL)}
  68(1):112--133.

\bibitem[{Denton \protect\BIBand{} Gupta(2003)}]{Denton2003}
Denton B, Gupta D (2003) A sequential bounding approach for optimal appointment
  scheduling. \emph{IIE Transactions} 35(11):1003--1016.

\bibitem[{Denton et~al.(2007)Denton, Viapiano, \protect\BIBand{}
  Vogl}]{Denton2007}
Denton B, Viapiano J, Vogl A (2007) Optimization of surgery sequencing and
  scheduling decisions under uncertainty. \emph{Health Care Management Science}
  10(1):13--24.

\bibitem[{El-Boghdadly et~al.(2020)El-Boghdadly, Nair, Pawa, \protect\BIBand{}
  Onwochei}]{El-Boghdadly2020}
El-Boghdadly K, Nair G, Pawa A, Onwochei DN (2020) Impact of parallel
  processing of regional anesthesia with block rooms on resource utilization
  and clinical outcomes: a systematic review and meta-analysis. \emph{Regional
  Anesthesia and Pain Medicine} 45(9):720--726.

\bibitem[{Erdogan \protect\BIBand{} Denton(2013)}]{Erdogan2013}
Erdogan S, Denton B (2013) Dynamic appointment scheduling of a stochastic
  server with uncertain demand. \emph{INFORMS Journal on Computing}
  25(1):116--132.

\bibitem[{Friedman et~al.(2006)Friedman, Sokal, Chang, \protect\BIBand{}
  Berger}]{Friedman2006}
Friedman DM, Sokal SM, Chang Y, Berger DL (2006) Increasing operating room
  efficiency through parallel processing. \emph{Annals of Surgery} 243:10--14.

\bibitem[{Gon{\c{c}}alves et~al.(2012)Gon{\c{c}}alves, Finardi,
  \protect\BIBand{} da~Silva}]{Goncalves2012}
Gon{\c{c}}alves RE, Finardi EC, da~Silva EL (2012) Applying different
  decomposition schemes using the progressive hedging algorithm to the
  operation planning problem of a hydrothermal system. \emph{Electric Power
  Systems Research} 83(1):19--27.

\bibitem[{Gul(2018)}]{Gul2018}
Gul S (2018) A stochastic programming approach for appointment scheduling under
  limited availability of surgery turnover teams. \emph{Service Science}
  10(3):277--288.

\bibitem[{Gul et~al.(2015)Gul, Denton, \protect\BIBand{} Fowler}]{Gul2015}
Gul S, Denton BT, Fowler JW (2015) A progressive hedging approach for surgery
  planning under uncertainty. \emph{INFORMS Journal on Computing}
  27(4):755--772.

\bibitem[{Gul et~al.(2011)Gul, Denton, Fowler, \protect\BIBand{}
  Huschka}]{Gul2011}
Gul S, Denton BT, Fowler JW, Huschka T (2011) Bi-criteria scheduling of
  surgical services for an outpatient procedure center. \emph{Production and
  Operations Management} 20(3):406--417.

\bibitem[{Gupta \protect\BIBand{} Denton(2008)}]{Gupta2008}
Gupta D, Denton B (2008) Appointment scheduling in health care: Challenges and
  opportunities. \emph{IIE transactions} 40(9):800--819.

\bibitem[{Hanss et~al.(2005)Hanss, Buttgereit, Tonner, Bein, Schleppers,
  Steinfath, Scholz, \protect\BIBand{} Bauer}]{Hanss2005}
Hanss R, Buttgereit B, Tonner P, Bein B, Schleppers A, Steinfath M, Scholz J,
  Bauer M (2005) Overlapping induction of anesthesia. \emph{Anesthesiology}
  103:391--400.

\bibitem[{Harders et~al.(2006)Harders, Melangoni, Weight, \protect\BIBand{}
  Sidhu}]{Harders2006}
Harders M, Melangoni M, Weight S, Sidhu T (2006) Improving operating room
  efficiency through process redesign. \emph{Anesthesiology} 140(4):509--516.

\bibitem[{Helseth(2016)}]{Helseth2016}
Helseth A (2016) Stochastic network constrained hydro-thermal scheduling using
  a linearized progressive hedging algorithm. \emph{Energy Systems}
  7(4):585--600.

\bibitem[{Hvattum \protect\BIBand{} L{\o}kketangen(2009)}]{Hvattum2009}
Hvattum LM, L{\o}kketangen A (2009) Using scenario trees and progressive
  hedging for stochastic inventory routing problems. \emph{Journal of
  Heuristics} 15(6):527.

\bibitem[{Khaniyev et~al.(2020)Khaniyev, Kayis, \protect\BIBand{}
  Gullu}]{Khaniyev2020}
Khaniyev T, Kayis E, Gullu R (2020) Next-day operating room scheduling with
  uncertain surgery durations: Exact analysis and heuristics. \emph{European
  Journal of Operational Research} 286(1):49--62.

\bibitem[{Laporte \protect\BIBand{} Louveaux(1993)}]{laporte1993integer}
Laporte G, Louveaux FV (1993) The integer l-shaped method for stochastic
  integer programs with complete recourse. \emph{Operations Research Letters}
  13(3):133--142.

\bibitem[{Lee \protect\BIBand{} Yih(2014)}]{Lee2014}
Lee S, Yih Y (2014) Reducing patient-flow delays in surgical suites through
  determining start-times of surgical cases. \emph{European Journal of
  Operational Research} 238:620--629.

\bibitem[{Mancilla \protect\BIBand{} Storer(2012)}]{Mancilla2012}
Mancilla C, Storer R (2012) A sample average approximation approach to
  stochastic appointment sequencing and scheduling. \emph{IIE Transactions}
  44(8):655--670.

\bibitem[{Mancilla \protect\BIBand{} Storer(2013)}]{Mancilla2013}
Mancilla C, Storer R (2013) Stochastic sequencing of surgeries for a single
  surgeon operating in parallel operating rooms. \emph{IIE Transactions on
  Healthcare Systems Engineering} 3(2):127--138.

\bibitem[{Marjamaa et~al.(2009)Marjamaa, Torkki, Hirvensalo, \protect\BIBand{}
  Kirvela}]{Marjamaa2009}
Marjamaa R, Torkki P, Hirvensalo E, Kirvela O (2009) What is the best workflow
  for an operating room? a simulation study of five scenarios. \emph{Health
  Care Management Science} 12:142--146.

\bibitem[{Rockafellar \protect\BIBand{} Wets(1991)}]{Rockafellar1991}
Rockafellar RT, Wets RJB (1991) Scenarios and policy aggregation in
  optimization under uncertainty. \emph{Mathematics of Operations Research}
  16(1):119--147.

\bibitem[{Sen \protect\BIBand{} Higle(2005)}]{sen2005c3}
Sen S, Higle JL (2005) The c3 theorem and a d2 algorithm for large scale
  stochastic mixed-integer programming: set convexification. \emph{Mathematical
  Programming} 104(1):1--20.

\bibitem[{Smith et~al.(2006)Smith, Sandberg, Foss, Massoli, Kanda, Barsoum,
  \protect\BIBand{} Schubert}]{Smith2008}
Smith MP, Sandberg WS, Foss J, Massoli K, Kanda M, Barsoum W, Schubert A (2006)
  High-throughput operating room system for joint arthroplasties durably
  outperforms routine processes. \emph{Anesthesiology} 109:25--35.

\bibitem[{Torkki et~al.(2005)Torkki, Marjamaa, Torkki, Kallio,
  \protect\BIBand{} Kirvela}]{Torkki2005}
Torkki P, Marjamaa R, Torkki M, Kallio P, Kirvela O (2005) Use of anesthesia
  induction rooms can increase the number of urgent orthopedic cases completed
  within 7 hours. \emph{Anesthesiology} 103:401--405.

\bibitem[{Van~Slyke \protect\BIBand{} Wets(1969)}]{VanSlykeWets1969}
Van~Slyke R, Wets R (1969) L-shaped linear programs with applications to
  optimal control and stochastic programming. \emph{SIAM Journal on Applied
  Mathematics} 17(4):638--663.

\bibitem[{Vandenberghe et~al.(2019)Vandenberghe, Vuyst, Aghezzaf,
  \protect\BIBand{} Bruneel}]{Vandenberghe2019}
Vandenberghe M, Vuyst SD, Aghezzaf EH, Bruneel H (2019) Surgery sequencing to
  minimize the expected maximum waiting time of emergent patients.
  \emph{European Journal of Operational Research} 275(3):971--982.

\bibitem[{Varmazyar et~al.(2020)Varmazyar, Akhavan-Tabatabaei, Salmasi,
  \protect\BIBand{} Modarres}]{Varmazyar2020}
Varmazyar M, Akhavan-Tabatabaei R, Salmasi N, Modarres M (2020) Operating room
  scheduling problem under uncertainty: Application of continuous phase-type
  distributions. \emph{IISE Transactions} 52(2):216--235.

\bibitem[{Watson \protect\BIBand{} Woodruff(2011)}]{Watson2011}
Watson JP, Woodruff DL (2011) Progressive hedging innovations for a class of
  stochastic mixed-integer resource allocation problems. \emph{Computational
  Management Science} 8(4):355--370.

\bibitem[{Yuan \protect\BIBand{} Sen(2009)}]{yuan2009enhanced}
Yuan Y, Sen S (2009) Enhanced cut generation methods for decomposition-based
  branch and cut for two-stage stochastic mixed-integer programs. \emph{INFORMS
  Journal on Computing} 21(3):480--487.

\bibitem[{Zhang \protect\BIBand{} Xie(2015)}]{Zhang2015}
Zhang Z, Xie X (2015) Simulation-based optimization for surgery appointment
  scheduling of multiple operating rooms. \emph{IIE Transactions}
  47(9):998--1012.

\end{thebibliography}

\appendix

\section{Proof of Proposition 1}

Letting $\mathbf{U} = \{u_{st}: s,t \in S\}$ denote the set of precedence relations from the first-stage of the two-stage SMIP model, suppose that $\sigma^F$ is the set of induction room assignments based on the first available induction room rule in Proposition 1. Furthermore, let $\sigma^* \neq \sigma^F$ denote the set of assignments in an optimal solution to the model defined by (1)-(19). The claim in Proposition 1 is that $\sigma^*$ can be converted to $\sigma^F$ without any change in the objective value.

Since $\sigma^* \neq \sigma^F$, there has to be at least one surgery patient in $\sigma^*$ for whom the induction room assignment will be different from that in $\sigma^F$. Let $i \in S$ be the earliest arriving such patient. Suppose that IR $n$ is the first one that becomes available after (or upon) the arrival of $i \in S$. By construction, this patient is assigned to a different IR, say $m$, which becomes available later. Let $k \in S$ be the first patient assigned to IR $n$ after the arrival of patient $i$, and patient $j \in S$ be the patient immediately preceding $k$ on the set of assignments to IR $n$.

\begin{figure}
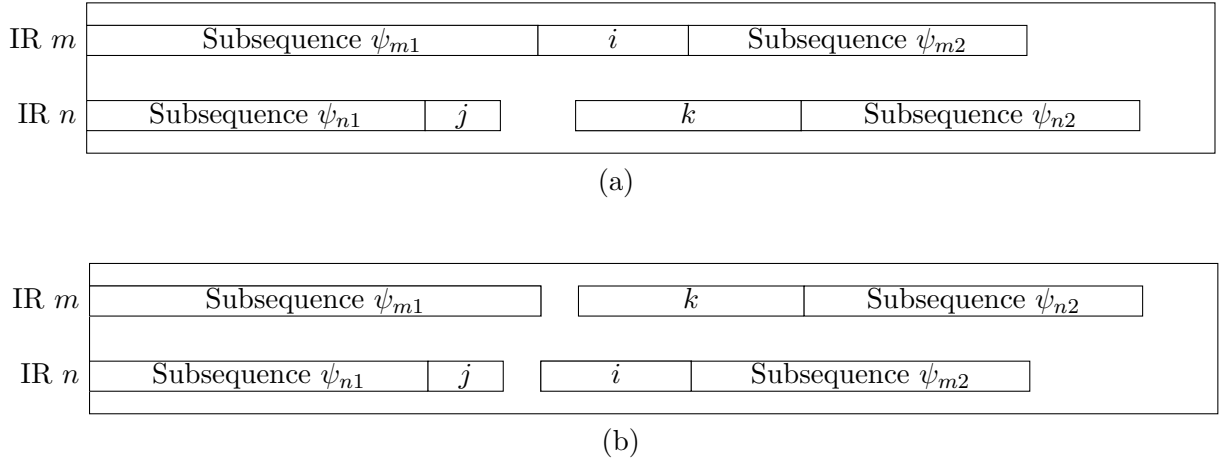

\begin{ganttchart}[inline,]{1}{30}
\ganttbar[inline=false]{IR $m$}{1}{12}
\ganttbar{Subsequence $\psi_{m1}$}{1}{12}
\ganttbar{$i$}{13}{16}
\ganttbar{Subsequence $\psi_{m2}$}{17}{25}
\\
\ganttbar{Subsequence $\psi_{n1}$}{1}{9}
\ganttbar{$j$}{10}{11}
\ganttbar[inline=false]{IR $n$}{14}{19}
\ganttbar{$k$}{14}{19}
\ganttbar{Subsequence $\psi_{n2}$}{20}{28}
\end{ganttchart}

\centering
(a)
\vspace{24pt}

\hspace{-0.4cm}
\begin{ganttchart}[inline,]{1}{30}
\ganttbar[inline=false]{IR $m$}{1}{12}
\ganttbar{Subsequence $\psi_{m1}$}{1}{12}
\ganttbar{$k$}{14}{19}
\ganttbar{Subsequence $\psi_{n2}$}{20}{28}
\\
\ganttbar{Subsequence $\psi_{n1}$}{1}{9}
\ganttbar{$j$}{10}{11}
\ganttbar[inline=false]{IR $n$}{13}{16}
\ganttbar{$i$}{13}{16}
\ganttbar{Subsequence $\psi_{m2}$}{17}{25}
\end{ganttchart}

\centering
(b)
\caption{The patient assignments for induction rooms IR $m$ and IR $n$ (a) under $\sigma^*$ and (b) after swapping}
\label{fig:proofchart}
\end{figure}

Figure \ref{fig:proofchart} illustrates this situation. Let Subsequences $\psi_{m1}$ and $\psi_{m2}$ represent the set of assignments (which may include idle times) to IR $m$ before and after the induction of patient $i$, respectively. Subsequences $\psi_{n1}$ and $\psi_{n2}$ denote the set of assignments to IR $n$ before the induction of patient $j$ and after that of patient $k$, respectively.

The induction of patient $j$ must finish before the start of induction for patient $i$ in IR $m$ in Figure \ref{fig:proofchart}, as otherwise this would be a contradiction to IR $n$ being the first available IR for patient $i$. Note further that by Constraints (8), we have $a_i < a_k$.

Consider the swapping of the assignments of patient $i$ and Subsequence $\psi_{m2}$ on IR $m$ with those of patient $k$ and Subsequence $\psi_{n2}$ on IR $n$. In this case:
\begin{itemize}
    \item Total IR idle time does not change, as the IR closure times have been swapped between these two IRs and total induction duration is the same.
    \item Total OR idle time is the same, as the OR schedules are unaffected (all patients leave their IR at the same times as in $\sigma^*$).
    \item Total patient waiting time is also identical, as the patients start their inductions at the same times as in $\sigma^*$.
\end{itemize}

To convert $\sigma^*$ to $\sigma^F$, one simply needs to find the earliest arriving patient whose IR assignment is different from $\sigma^F$, then find the first patient assigned to the IR that was first available for the former patient, and swap the subsequences as in Figure \ref{fig:proofchart} without changing the objective value. The process continues with the next earliest arriving patient with a different assignment, and so on. In the worst case, the conversion to $\sigma^F$ can be made in at most $|S|$ steps without any change in the objective. \QEDB

\end{document}